\documentclass[11pt]{article}
\usepackage{paper,asb}

\usepackage{cite}
\usepackage{soul}
\usepackage{bbm}
\usepackage{mathtools}

\newcommand{\papertitle}{Negative Curvature Methods with High-Probability
Complexity Guarantees for Stochastic Nonconvex
Optimization}
\newcommand{\paperauthora}{Albert S. Berahas}
\newcommand{\paperauthoraaffiliation}{Department of Industrial and Operations Engineering, University of Michigan}
\newcommand{\paperauthoraemail}{albertberahas@gmail.com}
\newcommand{\paperauthorb}{Wanping Dong}
\newcommand{\paperauthorbemail}{wanpingd@umich.edu}
\newcommand{\paperauthorc}{Raghu Bollapragada}
\newcommand{\paperauthorcaffiliation}{Operations Research and Industrial Engineering Program, University of Texas at Austin}
\newcommand{\paperauthorcemail}{raghu.bollapragada@utexas.edu}

\begin{document}
\title{\papertitle}

\author{\paperauthora \footnotemark[2] \footnotemark[1]
   \and \paperauthorc\footnotemark[3]\and \paperauthorb\footnotemark[2]\ }

\maketitle

\renewcommand{\thefootnote}{\fnsymbol{footnote}}
\footnotetext[2]{\paperauthoraaffiliation. (\url{\paperauthoraemail},\url{\paperauthorbemail})}
\footnotetext[3]{\paperauthorcaffiliation. (\url{\paperauthorcemail})}
\footnotetext[1]{Corresponding author.}
\renewcommand{\thefootnote}{\arabic{footnote}}

\begin{abstract}{
This paper develops negative curvature methods for continuous nonlinear unconstrained optimization in stochastic settings, in which function, gradient, and Hessian information is available only through probabilistic oracles, i.e., oracles that return approximations of a certain accuracy and reliability. We introduce conditions on these oracles and design a two-step framework that systematically combines gradient and negative curvature steps. The framework employs an early-stopping mechanism to guarantee sufficient progress and uses an adaptive mechanism based on an Armijo-type criterion to select the step sizes for both steps. We establish high-probability iteration-complexity guarantees for attaining second-order stationary points, deriving explicit tail bounds that quantify the convergence neighborhood and its dependence on oracle noise. Importantly, these bounds match deterministic rates up to noise-dependent terms, and the framework recovers the deterministic results as a special case. Finally, numerical experiments demonstrate the practical benefits of exploiting negative curvature directions even in the presence of noise.

}
\end{abstract}

\section{Introduction}

In this paper, we consider unconstrained nonlinear optimization problems of the form, $\min_{x \in \mathbb{R}^{n}} f(x)$, 
where $f:\mathbb{R}^{n} \rightarrow \mathbb{R}$ is twice continuously differentiable and potentially nonconvex. We study settings in which the exact objective function and its associated derivatives ($\nabla f$ and $\nabla^2 f$) are not available; instead, zeroth-, first-, and second-order approximations are obtained through probabilistic oracles, i.e., oracles that return approximations of a certain accuracy and reliability. Such problems 
arise 
in applications such as simulation optimization \cite{pasupathy2018sampling}, machine learning  \cite{bottou2018optimization}, and decision-making \cite{rockafellar2000optimization}.

While most of the literature in the aforementioned probabilistic-oracle setting focuses on algorithms that guarantee convergence to (approximate) first-order stationary points, our goal is to develop, analyze, and implement methods capable of converging, at appropriate rates, to (approximate) second-order stationary points, defined as points $x \in \mathbb{R}^n$ and constants $\bar{\epsilon}_g,\bar{\epsilon}_{\lambda},\bar{\epsilon}_g>0$ such that 
\begin{align}\label{eq.2ndstat}
    \|\nabla f(x)\|_2 < \bar{\epsilon}_g, \quad \text{and} \quad  \lambda_{\min}(\nabla f^2 (x)) > -\max\{\bar{\epsilon}_{\lambda}, \bar{\epsilon}_H\}.
\end{align}
To this end, we develop adaptive stochastic algorithms that exploit negative curvature directions and operate solely using the outputs of probabilistic zeroth-, first-, and second-order oracles. The proposed methods follow a two-step structure and alternate between gradient (or general descent) steps and negative curvature steps \cite{berahas2024exploiting, curtis2019exploiting}, employ an Armijo-type step-search mechanism to adaptively select step sizes under noisy evaluations \cite{berahas2021global, berahas2025sequential, jin2021high}, and are endowed with high-probability iteration-complexity guarantees. In particular, we show that, after a sufficient number of iterations, the algorithms reach an $(\bar{\epsilon}_g, \bar{\epsilon}_{\lambda}, \bar{\epsilon}_H)$-stationary point \eqref{eq.2ndstat} with overwhelmingly high probability, where the resulting convergence neighborhood scales with the oracle noise and matches the deterministic setting up to noise-dependent terms.

\subsection{Related Literature}

In this section, we review the literature most closely related to our problem setting and the components of our algorithmic framework. We focus on prior work on negative curvature methods, optimization under various inexact or stochastic oracle models, and stochastic line/step-search strategies. These areas form the foundation on which our proposed method is built.

A large body of work has established the role of negative curvature (NC) in escaping saddle points and achieving second-order convergence in deterministic nonconvex optimization \cite{nocedal1999numerical}. Classical methods include curvilinear and two-step schemes that combine descent and NC directions \cite{forsgren1995computing, goldfarb1980curvilinear, more1979use}, as well as recent variants that exploit NC directions explicitly at each iteration \cite{berahas2024exploiting, curtis2019exploiting}. Trust-region and cubic-regularization methods also implicitly exploit NC information through their subproblem solvers (e.g., truncated CG or Lanczos) and enjoy well-understood global convergence and complexity guarantees \cite{cartis2011adaptive, nesterov2006cubic}. These deterministic developments motivate the structure of our algorithm but rely on access to exact or deterministically bounded function and derivative information, assumptions that break down in stochastic environments.

NC methods have been studied in noisy settings, e.g.,  deterministic inexact oracles with bounded errors \cite{berahas2024exploiting,li2025randomized, royer2018complexity}, and stochastic settings in which convergence is established in expectation \cite{allen2018natasha, berahas2024exploiting, curtis2019exploiting}. However, most of these approaches require restrictive assumptions on accuracy and typically guarantee only first-order convergence or convergence in expectation rather than with high probability. To the best of our knowledge, no prior work provides high-probability second-order complexity guarantees for NC-based methods under general probabilistic error oracle settings. 

The study of optimization methods in the presence of inexact oracles has evolved to include deterministic \cite{berahas2019derivative, carter1991global}, in expectation \cite{berahas2024exploiting, bottou2018optimization, paquette2020stochastic, royer2018complexity}, and probabilistic \cite{bellavia2022linesearch, berahas2025sequential, cao2023first,jin2021high} error models. Probabilistic oracles, in particular, allow the returned estimates to be inaccurate or biased with nontrivial probability and are increasingly used to model realistic stochastic evaluations. However, existing works primarily focus on zeroth- and first-order information, and probabilistic second-order oracles suitable for reliably identifying NC directions in noisy settings have received little attention. Our setting aligns with this general probabilistic framework but extends it by considering both bounded noise and subexponential noise models for function evaluations, biased gradient approximations that may be arbitrarily inaccurate with nonzero probability, and probabilistic Hessian-vector information for detecting NC.

Stochastic line-search and step-search methods provide another important ingredient for our work. Traditional deterministic line-search procedures rely critically on accurate function and derivative information and are not directly applicable in noisy settings \cite{nocedal1999numerical}. To address this, several adaptive schemes have been proposed, including stochastic backtracking strategies with sample-size adjustment \cite{bollapragada2018adaptive,friedlander2012hybrid}. More recently, step-search methods, also called stochastic line-search methods, have been developed to adaptively increase or decrease step sizes based on probabilistic sufficient-decrease conditions, re-evaluating the oracle at the same iterate when a trial step fails; see e.g.,  \cite{berahas2021global, cartis2018global, jin2021high, paquette2020stochastic}. These strategies yield improved robustness in stochastic environments but largely target first-order methods and do not address the handling of NC directions or provide second-order convergence guarantees.

In summary, while negative curvature methods, inexact or probabilistic oracles, and stochastic step-search techniques have each been studied in isolation, there remains a gap in combining these components to obtain second-order guarantees in general noisy settings. Our work fills this gap by developing a negative-curvature-based adaptive step-search framework tailored to general probabilistic settings and endowed with high-probability iteration-complexity guarantees.

\subsection{Contributions} 

In this work, we develop a negative curvature method for stochastic nonconvex optimization under a general probabilistic-oracle setting, in which zeroth-, first-, and second-order information is available only through noisy approximations of prescribed accuracy and reliability. The framework accommodates a wide range of noise models, function evaluations may contain either bounded noise or noise with subexponential tails, gradient estimates may be biased and arbitrarily inaccurate with controlled probability, and Hessian information is obtained through probabilistic second-order oracles tailored to detecting negative curvature under noise.

Our first contribution is a flexible algorithmic framework for this setting. The framework follows a two-step structure alternating between descent (gradient-related) and negative curvature directions, and relaxes traditional sufficient-decrease conditions to account for noise. An Armijo-type step-search procedure is employed for both step types, enabling adaptive step-size selection and re-evaluation of search directions when needed. The framework also incorporates an efficient mechanism for selecting negative curvature directions using only two function evaluations and no gradient evaluations, reducing computational cost and enhancing robustness in large-scale settings. Overall, the scheme encompasses a broad class of stochastic optimization methods as special cases while introducing new components essential for handling noise.

Our second contribution is a high-probability second-order convergence and complexity analysis for the proposed method. We establish explicit tail bounds showing that the probability of requiring more than $\mathcal{O}(\max\{\bar{\epsilon}_g^{-2}, \bar{\epsilon}_H^{-3}, \bar{\epsilon}_{\lambda}^{-3}\})$ iterations to reach an ($\bar{\epsilon}_g$, $\bar{\epsilon}_H$, $\bar{\epsilon}_{\lambda}$)-second-order stationary point decays exponentially with the number of iterations. These guarantees match the deterministic rates up to noise-dependent terms, recover deterministic results as a special case when noise vanishes, and extend prior stochastic analyses by operating under broader probabilistic oracle assumptions. 

Finally, we implement a practical method motivated by our theoretical framework, replacing idealized conditions and second-order computations with an efficient conjugate-gradient–based subsolver to detect negative curvature, while retaining the step-selection and acceptance mechanisms motivated by the analysis. Preliminary experiments on a classical test problem illustrate the robustness of the approach to noise, often outperforming baseline stochastic methods in noisy settings.

\subsection{Paper Organization}

The paper is organized as follows. We conclude this section by introducing the notation used throughout the paper. 
The assumptions, probabilistic oracles, and the algorithm we propose are presented in Section \ref{sec.assump_alg}. In Section~\ref{sec.analysis}, we establish high-probability complexity bounds for reaching second-order stationary points. Section~\ref{sec.experiment} reports numerical experiments that illustrate the practical performance of the algorithm. Finally, concluding remarks are provided in Section~\ref{sec.remarks}.

\subsection{Notation}

Let $\mathbb{N}$ and $\mathbb{R}$ denote the sets of natural and real numbers, respectively. For $n,m \in \mathbb{N}$, let $\mathbb{R}^n$ and $\mathbb{R}^{n\times m}$ denote the sets of $n$-dimensional real vectors and $n\times m$ real matrices, respectively. For a matrix $A$, let $\lambda_{\min}(A)$ ($\lambda_{\max}(A)$) denote its minimum (maximum) eigenvalues. We consider a sequence of iterates $\{x_k, \hat{x}_k\}_{k \in \mathbb{N}}$ with $x_k, \hat{x}_k \in \mathbb{R}^n$. At iteration $k$, the algorithm may generate a descent direction $d_k \in \mathbb{R}^n$ and a negative curvature direction $p_k \in \mathbb{R}^n$, with corresponding step sizes $\alpha_k > 0$ and $\beta_k > 0$, respectively. For each $k \in \mathbb{N}$, let $F_k$ denote an estimate of $f(x_k)$, $g_k \in \mathbb{R}^n$ an estimate of $\nabla f(x_k)$, and $H_k \in \mathbb{R}^{n \times n}$ an estimate of $\nabla^2 f(\hat{x}_k)$. 

\section{Assumptions, Oracles, and Algorithm} \label{sec.assump_alg}
In this section, we present the assumptions, define the oracles, and describe the algorithmic framework. Throughout the paper, we make the following assumption about the objective function. This assumption is standard in the literature of methods that exploit negative curvature, and more generally in the analysis of second-order methods; see e.g., \cite{cao2023first, curtis2019exploiting, li2025randomized}.
\bassumption \label{asm.f.general}
    The function $f: \mathbb{R}^{n} \to \mathbb{R}$ is continuously differentiable, the gradient of $f$ is $L_g$-Lipschitz continuous, and the Hessian of $f$ is $L_H$-Lipschitz continuous, for all $x \in \mathbb{R}^{n}$. Moreover, the function $f$ is bounded below by a scalar $\bar{f} \in \mathbb{R}$.
\eassumption

Our algorithm relies on approximations of the objective function and its derivatives provided by probabilistic oracles satisfying prescribed accuracy and reliability conditions. Precise definitions of the zeroth-, first-, and second-order oracles follow.


\boracle \label{def.0th-order_orc}
    Given a point $x \in \mathbb{R}^{n}$, the oracle computes $F(x, \xi^{(0)})$, a (random) estimate of the function value $f(x)$, where $\xi^{(0)}$ is a random variable defined on a probability space $(\Xi, \mathcal{F}, \P)$ whose distribution may depend on $x$. 
    For any $x \in \mathbb{R}^{n}$, $e(x, \xi^{(0)}) = |F(x, \xi^{(0)}) - f(x)|$ satisfies at least one of the two conditions:
    \begin{enumerate}[leftmargin=0.5cm]
        \item \label{orc.0th-order_bounded}(Deterministically bounded noise) There is a constant $\epsilon_f \geq 0$ such that $e(x, \xi^{(0)}) \leq \epsilon_f$ for all realizations of $\xi^{(0)}$.
        \item \label{orc.0th-order_subexp_orc} (Independent subexponential noise) There are constants $\epsilon_f \geq 0$ and $a > 0$ such that for all $s > 0$, $\P_{\xi^{(0)}} \left[e(x, \xi^{(0)}) \geq s\right] \leq \exp(-a(s-\epsilon_f))$.
    \end{enumerate}
\eoracle
\bremark
When Oracle~\ref{def.0th-order_orc}.\ref{orc.0th-order_bounded} is queried at a point $x \in \mathbb{R}^n$, an estimate $F(x,\xi^{(0)})$ is returned whose error is deterministically bounded by $\epsilon_f$. Under Oracle~\ref{def.0th-order_orc}.\ref{orc.0th-order_subexp_orc}, the function noise has subexponential tails with a bias allowance of $\epsilon_f$. 
The parameters $a > 0$ and $\epsilon_f \ge 0$ are global constants independent of $x$. Similar formulations of zeroth-order probabilistic oracles can be found in \cite{cao2023first,jin2021high}.
\eremark

\boracle \label{def.1st-order_orc}
Given a probability $p_g \in \left(\tfrac{1}{2}, 1\right]$ and a point $x \in \mathbb{R}^{n}$, the oracle computes $g(x, \xi^{(1)})$ a (random) estimate of the gradient $\nabla f(x)$ that satisfies
\begin{equation*} 
    \P_{\xi^{(1)}} \left[\|g(x, \xi^{(1)}) - \nabla f(x)\|_2 \leq \epsilon_g + \kappa_g \|\nabla f(x)\|_2 \right] \geq p_g
\end{equation*}
where $\epsilon_g \geq 0$,  $\kappa_g \geq 0$, and $\xi^{(1)}$ is a random variable defined on a probability space $(\Xi, \mathcal{F}, \P)$ whose distribution may depend on $x$. 
\eoracle
\bremark
Oracle~\ref{def.1st-order_orc} returns a gradient estimate $g(x,\xi^{(1)})$ satisfying the accuracy condition with probability at least $p_g > \tfrac{1}{2}$. The right-hand side consists of absolute and relative terms, where $\epsilon_g$ and $\kappa_g$ are nonnegative constants intrinsic to the oracle. When $\|\nabla f(x)\|_2$ is large or $\epsilon_g = 0$, the relative term dominates and the condition enforces norm-relative accuracy \cite{bollapragada2018adaptive, byrd2012sample, carter1991global}. When $\|\nabla f(x)\|_2$ is small or $\kappa_g = 0$, the absolute term dominates. The constant $\epsilon_g$ bounds the best achievable gradient accuracy (see \eqref{eq.value.eps}). If $p_g = 1$, the condition is deterministic. This combined absolute/relative probabilistic accuracy condition is standard; see, e.g., \cite{berahas2024exploiting, berahas2025sequential, bollapragada2018adaptive, jin2021high, li2025randomized}.
\eremark

Prior to defining the second-order probabilistic oracle, we define a negative curvature direction $q \in \mathbb{R}^n$ for a matrix $H \in \mathbb{R}^{n\times n}$.
\bdefinition \label{cond.q}
Consider a matrix $H \in \mathbb{R}^{n\times n}$ for which $\lambda_{\min}(H) < 0$. A vector $q \in \mathbb{R}^n$ is a direction of negative curvature  of $H$ if 
\begin{subequations}\label{fom.q}
    \begin{align}
         q^{\mathsf{T}} H q & \leq \gamma \lambda_{\min}(H) \|q\|_2^2 < 0, \quad \gamma \in (0, 1] \label{fom.q1} \\
         \|q\|_2 &= \delta |\lambda_{\min}(H)|, \quad\quad \quad \quad\delta \in (0, \infty).  \label{fom.q2}
    \end{align}
\end{subequations}
\edefinition
The second-order probabilistic oracle can then be defined as follows.
\boracle \label{def.2nd-order_orc}
Given a probability $p_H \in (0.5, 1]$ and a point $x \in \mathbb{R}^{n}$, the oracle computes $H(x, \xi^{(2)})$, a (random) estimate of the Hessian $\nabla^2 f(x)$, such that
\begin{equation*} 
    \begin{split}
        & \P_{\xi^{(2)}}\left[\|(\nabla^2 f(x) - H(x, \xi^{(2)})) q\|_2 \leq \epsilon_H + \kappa_H |\lambda_{\min}(H(x, \xi^{(2)}))|\|q\|_2, \right. \\
        & \qquad\qquad  \left. | \lambda_{\min}(\nabla^2 f(x)) - \lambda_{\min}(H(x, \xi^{(2)}))| \leq \epsilon_{\lambda} + \kappa_{\lambda} |\lambda_{\min}(\nabla^2 f(x))| \right] \geq p_H.
    \end{split}
\end{equation*}
where $\epsilon_H \geq 0$,  $\kappa_H \geq 0$, $\epsilon_{\lambda} \geq 0$,  $\kappa_{\lambda} \geq 0$, $\xi^{(2)}$ is a random variable defined on a probability space $(\Xi, \mathcal{F}, \P)$ whose distribution may depend on the input $x$, and $q$ is a negative curvature direction of $H(x, \xi^{(2)})$.
\eoracle
\bremark
Oracle~\ref{def.2nd-order_orc} imposes two inexactness conditions, enforced only when negative curvature is detected. As in Oracle~\ref{def.1st-order_orc}, each combines relative and absolute terms with nonnegative constants intrinsic to the oracle. The first ensures directional accuracy, requiring the Hessian estimate to be accurate only along a negative curvature direction \cite{berahas2024exploiting, xu2020newton}. The second controls the approximation of the left-most eigenvalue. Together, these conditions are strictly weaker than uniform spectral-norm bounds, i.e., $\|\nabla^2 f(x) - H(x,\xi^{(2)})\|_2$. If $p_H = 1$, they reduce to the inexact deterministic setting.
\eremark

For notational simplicity, we denote the realizations of random estimates at each iterate $x_k$ generated by the algorithm by $F_k = F (x_k, \xi_k^{(0)})$, $F_k^+ = F(x_k + \alpha_k d_k, \xi_k^{(0+)})$, $\hat{F}_k = F (\hat{x}_k, \hat{\xi}_k^{(0)})$, $\hat{F}_k^{\pm} = F(\hat{x}_k {\pm} \beta_k q_k, \hat{\xi}_k^{(0{\pm})})$, $g_k = g(x_k, \xi_k^{(1)})$, $H_k = H(\hat{x}_k, \hat{\xi}_k^{(2)})$, and $\lambda_k = \lambda_{\min} (H_k)$. The algorithmic framework is given in Algorithm \ref{alg.theory}.
\begin{algorithm}[]
    \caption{Step Search Negative Curvature Method with Probabilistic Oracles} \label{alg.theory}
    \begin{algorithmic}[1]
        \Require{$x_{0} \in \mathbb{R}^{n}$, $\alpha_0, \beta_0 > 0$, $\tau \in (0,1]$, $c_d, c_p \in (0, 1)$, $e_f, c_g, c_H \geq 0$, $\bar{\epsilon}_g, \bar{\epsilon}_H, \bar{\epsilon}_{\lambda} \geq 0$}
        \ForAll{$k \in \{0, 1, \dots\}\eqqcolon\N{}$} 
            \State generate $g_k$ via Oracle \ref{def.1st-order_orc} \label{line.start_dc}
            \If{$\|g_k\|_2 \leq c_g \bar{\epsilon}_g$} 
                \ set $\hat{x}_{k} \leftarrow x_k$, $\alpha_{k+1} \leftarrow \alpha_k$ \label{line.term_dc}
            \Else \ set $d_k \leftarrow - g_k$ and generate $F_k$ and $F_k^+$ via Oracle \ref{def.0th-order_orc}
                \If{$F_k^+ \leq F_k + c_d \alpha_k d_k^{\mathsf{T}} g_k + e_f$} \ set $\hat{x}_{k} \leftarrow x_k + \alpha_k d_k$, $\alpha_{k+1} \leftarrow \tau^{-1} \alpha_k$ \label{line.armijo_dc} 
                \Else\ set $\hat x_k \gets x_k$, $\alpha_{k+1} \gets \tau\alpha_k$ \label{line.ss2_dc}
                \EndIf
            \EndIf \label{line.end_dc}
            \State generate $H_k$ via Oracle \ref{def.2nd-order_orc} and compute $\lambda_k$ \label{line.start_nc}
            \If{$\lambda_k \geq - c_H \max\{\bar{\epsilon}_H, \bar{\epsilon}_{\lambda}\}$} \ set $x_{k+1} \leftarrow \hat{x}_k$, $\beta_{k+1} \leftarrow \beta_k$ \label{line.term_nc}
            \Else \ compute negative curvature direction $q_k$ satisfying Condition \ref{cond.q} 
            
            and generate $\hat{F}_k$, $\hat{F}_k^+$ and $\hat{F}_k^-$ via Oracle \ref{def.0th-order_orc}
                \If{$\min\{\hat{F}_k^+, \hat{F}_k^-\} \leq \hat{F}_k + c_p \beta_k^2 q_k^{\mathsf{T}} H_k q_k + e_f$} \label{line.armijo_nc}
                    \State set $p_k \leftarrow \arg\min_{ \omega \in \{-1,1\}}\{F(\hat{x}_k + \omega \beta_k q_k, \hat{\xi}_k^{(0\pm)})\} q_k$ \label{line.sign_selection}
                    \State set $x_{k+1} \leftarrow \hat{x}_k + \beta_k p_k$, $\beta_{k+1} \leftarrow \tau^{-1} \beta_{k}$ \label{line.ss1_nc}
                \Else \ set $x_{k+1} \leftarrow \hat{x}_k$, $\beta_{k+1} \leftarrow \tau \beta_{k}$ \label{line.ss2_nc}
                \EndIf
            \EndIf \label{line.end_nc}
	\EndFor
    \end{algorithmic}
\end{algorithm}

\vspace{-0.25cm}

\bremark
We make a few remarks about Algorithm~\ref{alg.theory}. 
\begin{itemize}[leftmargin=0.5cm]
    \item Algorithm~\ref{alg.theory} alternates between descent and negative curvature steps, with the latter computed only when negative curvature is detected. As presented $d_k \leftarrow -g_k$, though more general forms $d_k \leftarrow -W_k g_k$, where $w_1 I \preceq W_k \preceq w_2 I$ for some $0 < w_1 \leq w_2$, are allowed without affecting the theoretical guarantees.
    \item The cost of the descent step (Lines~\ref{line.start_dc}--\ref{line.end_dc}) is dominated by the computation of a gradient estimate $g_k$ and two function evaluations, $F_k$ and $F_k^+$. In contrast, the cost of a negative curvature step (Lines~\ref{line.start_nc}--\ref{line.end_nc}) arises from computing a Hessian estimate $H_k$, its minimum eigenvalue $\lambda_k$, and potentially a negative curvature direction $q_k$, together with three function evaluations, $\hat{F}_k$ and $\hat{F}_k^{\pm}$.

    \item The early-termination rules in Lines~\ref{line.term_dc} and~\ref{line.term_nc} accept only steps whose predicted decrease, based on the estimated gradient norm or magnitude of negative curvature, is sufficiently large, thereby filtering out steps unlikely to produce meaningful progress. These rules do not terminate the algorithm; when triggered, the corresponding descent or negative curvature step is skipped. Enforcing these conditions requires knowledge of, or reliable estimates for, the absolute error parameters of the first- and second-order oracles. Early termination is necessary because our framework involves two step types and a second-order stopping criterion with two conditions, requiring separate control of successful and large-step events in the stopping-time analysis, in contrast to prior works that consider either a stopping time defined by a single condition \cite{jin2021high} or a single step per iteration \cite{cao2023first}.
    \item Step sizes are determined via a step-search procedure and accepted using relaxed Armijo conditions (Lines~\ref{line.armijo_dc}--\ref{line.ss2_dc} and~\ref{line.armijo_nc}--\ref{line.ss2_nc}). This stochastic line search repeatedly re-evaluates the required oracle information at the current iterate until the condition is satisfied, improving robustness under noisy function evaluations. The relaxed conditions include a noise-tolerance parameter $e_f$ to ensure acceptable step sizes despite function noise. For negative curvature steps, the Armijo condition is further modified to account for curvature effects, as in \cite{goldfarb1980curvilinear, mccormick1977modification, more1979use}.
    \item A key adaptation concerns the selection of the negative curvature direction. Since both $q_k$ and $-q_k$ are valid negative curvature directions, the goal is to choose the sign that yields descent. Verifying descent typically requires gradient information, increasing the cost of a negative curvature step. Prior approaches either select the sign at random \cite{li2025randomized} or use approximate gradients \cite{berahas2024exploiting}, trading per-iteration cost against effective curvature exploitation. In contrast, our method determines the sign by comparing two trial function evaluations and selecting $\min\{\hat F_k^+, \hat F_k^-\}$ (Line~\ref{line.armijo_nc}), requiring at most one additional function evaluation. This provides an efficient and reliable sign-selection mechanism without gradient information and integrates naturally with the Armijo conditions in the step-search procedure.
\end{itemize}
\eremark

\vspace{-0.1cm}

To analyze the convergence of Algorithm \ref{alg.theory}, we 
introduce the following notation.
Let $\{X_k\}$, $\{X_k^+\}$, $\{\hat{X}_k\}$, $\{\hat{X}_k^+\}$, and $\{\hat{X}_k^-\}$ be random vectors in $\mathbb{R}^{n}$ with realizations $x_k$, $x_k + \alpha_k d_k$, $\hat{x}_k$, $\hat{x}_k + \beta_k q_k$, and $\hat{x}_k - \beta_k q_k$, respectively. We define the magnitudes of the zeroth-order oracle errors as random variables $\{E_k\}$, $\{E_k^+\}$, $\{\hat{E}_k\}$, $\{\hat{E}_k^+\}$, and $\{\hat{E}_k^-\}$ by $E_k \coloneqq | f(X_k) - F(X_k, \xi_k^{(0)})|$, $E_k^+ \coloneqq | f(X_k + \alpha_k d_k) - F(X_k + \alpha_k d_k, \xi_k^{(0+)})|$, $\hat{E}_k \coloneqq | f(\hat{X}_k) - F(\hat{X}_k, \hat{\xi}_k^{(0)})|$, $\hat{E}_k^+ \coloneqq | f(\hat{X}_k + \beta_k q_k) - F(\hat{X}_k + \beta_k q_k, \hat{\xi}_k^{(0+)})|$, $\hat{E}_k^- \coloneqq | f(\hat{X}_k - \beta_k q_k) - F(\hat{X}_k - \beta_k q_k, \hat{\xi}_k^{(0-)})|$, respectively, with realizations 
$e_k$, $e_k^+$, $\hat{e}_k$, $\hat{e}_k^+$, and $\hat{e}_k^-$. 
For notational convenience, we let $\{\alpha_k\}$ and $\{\beta_k\}$ denote the random step size sequences, with realizations written using the same symbols.

At iterate $X_k$ of Algorithm \ref{alg.theory}, randomness arises in the following operations:
\begin{enumerate}[leftmargin=0.5cm]
    \item \textbf{Stochastic gradient at $X_k$:} Given $X_k$ a gradient estimate $g_k = g(X_k, \xi_k^{(1)})$ is computed via Oracle \ref{def.1st-order_orc}. Then, given $X_k$, $\alpha_k$, and $g_k$, the intermediate point $X_k^+$ is set deterministically.
    \item \textbf{Stochastic function values at $X_k$ and $X_k^+$:} Given $X_k$ and $X_k^+$ function estimates are computed via Oracle \ref{def.0th-order_orc}, introducing the corresponding errors errors $E_k$ and $E_k^+$. These estimates are used to update $\alpha_{k+1}$ and form $\hat{X}_k$.
    \item \textbf{Stochastic Hessian at $\hat{X}_k$:} Given $\hat{X}_k$ a Hessian approximation $H(\hat{X}_k, \hat{\xi}_k^{(2)})$ is computed via Oracle \ref{def.2nd-order_orc}, and then the smallest eigenvalue $\lambda_k$ and an associated negative curvature direction $q_k$ (if one exists) are obtained deterministically. If $q_k$ exists, trial points $\hat{X}_k^+$ and $\hat{X}_k^-$ are constructed using $\hat{X}_k$, $q_k$ and $\beta_k$. Otherwise, the algorithm sets $X_{k+1} \gets \hat{X}_k$ and $\beta_{k+1} \gets \beta_{k}$. 
    \item \textbf{Stochastic function values at $\hat{X}_k$, $\hat{X}_k^+$ and $\hat{X}_k^-$:} When $q_k$ exists, given $\hat{X}_k$, $\hat{X}_k^+$ and $\hat{X}_k^-$ function estimates are computed via Oracle \ref{def.0th-order_orc}. These estimates are used to chose between the two directions and thereby determine $X_{k+1}$ and $\beta_{k+1}$.
\end{enumerate}
Thus, conditioned on $X_k, \alpha_k, \beta_k$, all randomness in iteration $k$ is captured by $\xi_k^{(0)}$, $\xi_k^{(0+)}$, $\xi_k^{(1)}$, $\hat{\xi}_k^{(0)}$, $\hat{\xi}_k^{(0\pm)}$, and $\hat{\xi}_k^{(2)}$. 
Let $\mathcal{Z}_{k} \coloneqq \left\{\xi_k^{(0)}, \xi_k^{(0+)}, \xi_k^{(1)}, \hat{\xi}_k^{(0)}, \hat{\xi}_k^{(0\pm)}, \hat{\xi}_k^{(2)}\right\}$ and $\mathcal{F}_{k}$ denote the $\sigma$-algebra representing the history of the algorithm up to iteration $k$, i.e., $\mathcal{F}_{k} \coloneqq \sigma \left( \cup_{i=1}^{k} \mathcal{Z}_i \right)$. 
For analysis purposes, we introduce the augmented filtrations, $\mathcal{F}_{k}' \coloneqq \sigma \left( \mathcal{F}_{k} \cup \sigma\left(\xi_{k+1}^{(1)} \right) \right)$ and $\hat{\mathcal{F}}_{k}' \coloneqq \sigma \left( \mathcal{F}_{k} \cup \sigma\left(\xi_{k+1}^{(0)}, \xi_{k+1}^{(0+)}, \xi_{k+1}^{(1)}, \xi_{k+1}^{(2)} \right) \right)$. 

Our analysis relies on categorizing iterations $k = 0, 1, \dots, t-1$, for $t\geq 1$, into different types. These types can be defined in terms of estimation accuracy using the following random indicator variables: 
\begin{align*}
    I_k^g & \coloneqq \mathbbm{1} \left\{ \|g(X_k, \xi_k^{(1)}) - \nabla f(X_k)\|_2 \leq \epsilon_g + \kappa_g \|\nabla f(X_k)\|_2 \right\}, \\
    I_k^H & \coloneqq \mathbbm{1} \left\{ \|(\nabla^2 f(X_k) - H(X_k, \xi_k^{(2)})) q_k\|_2 \leq \epsilon_H^2 + \kappa_H |\lambda_k|\|q_k\|_2, \right. \\
    & \left. \qquad \quad  |\lambda_{\min}(H(X_k, \xi_k^{(2)})) - \lambda_{\min}(\nabla^2 f(X_k))| \leq \epsilon_{\lambda} + \kappa_{\lambda} |\lambda_{\min}(\nabla^2 f(X_k))| \right\},\\
    \text{and} \quad I_k^f & \coloneqq \mathbbm{1} \left\{ E_k + E_k^+ \leq e_f \right\}, \quad \hat{I}_k^f \coloneqq \mathbbm{1} \left\{ \hat{E}_k + \max\{\hat{E}_k^+, \hat{E}_k^-\} \leq e_f \right\}.
\end{align*}
With high probability, the number of iterations with accurate function evaluations being generated has a lower bound with respect to $t$. This holds trivially for Oracle~\ref{def.0th-order_orc}.\ref{orc.0th-order_bounded}.
\blemma \label{lmm.f.accuracy}
    Let $p_f = 1 - 3 \exp(-a(\tfrac{e_f}{2} - \epsilon_f) )$, 
    where $a$ is the positive constant from Oracle \ref{def.0th-order_orc}.\ref{orc.0th-order_subexp_orc}. Then, the indicators variables $I_k^f$ and $\hat{I}_k^f$ satisfy the submartingale condition, i.e.,  $\P\left[I_k^f = 1 | \mathcal{F}_{k-1}' \right] \geq p_f$ and $\P \left[\hat{I}_k^f = 1 | \hat{\mathcal{F}}_{k-1}' \right] \geq p_f$. 
    Furthermore, if $e_f > 2 \epsilon_f + \tfrac{2}{a} \log 6$ (equivalently, $p_f > \tfrac{1}{2}$), then for any positive integer $t \geq 1$ and any $\bar{p}_f \in [0, p_f]$, it follows that $\P\left[ \sum_{k=0}^{t-1} I_k^f \geq \bar{p}_f t \right] \geq 1-\exp \left(- \tfrac{(\hat{p}_f - p_f)^2}{2p_f^2} t \right)$.
\elemma
\bproof
By the definition of Oracle \ref{def.0th-order_orc}.\ref{orc.0th-order_subexp_orc}, 
at the descent steps, for all $k \in \N{}$
\begin{equation} \label{eq.prob_bound_Ek}
    \P_{\xi_k^{(0)}}\left[ E_k > s | \mathcal{F}_{k-1}' \right] \leq \exp (a(\epsilon_f - s)), \quad \P_{\xi_k^{(0+)}}\left[ E_k^+ > s | \mathcal{F}_{k-1}' \right] \leq \exp (a(\epsilon_f - s)).
\end{equation}
Consequently, by the definition of $I_k^f$, Oracle \ref{def.0th-order_orc}.\ref{orc.0th-order_subexp_orc}, and $p_f$, it follows that
\begin{align*}
    \P \left[I_k^f = 0| \mathcal{F}_{k-1}' \right]  = \P \left[E_k + E_k^+ > e_f| \mathcal{F}_{k-1}' \right] 
    & \leq \P \left[E_k > e_f/2 \text{ or } E_k^+ > e_f/2 | \mathcal{F}_{k-1}'\right] \\
    & \leq \P \left[E_k > e_f/2 | \mathcal{F}_{k-1}'\right] + \P \left[ E_k^+ > e_f/2 | \mathcal{F}_{k-1}'\right] \\
    & \leq 2 \exp (a(\epsilon_f - e_f/2)) \leq 1- p_f.
\end{align*}
Thus, 
the random process $\left\{ \sum^{t-1}_{k=0} I_k^f - \bar{p}_f t \right\}_{t=0,1,\dots}$ is a submartingale.
Since
\begin{align*}
    \left| \left( \sum^{(t+1)-1}_{k=0} I_k^f - \bar{p}_f(t + 1) \right) - \left(\sum^{t-1}_{k=0} I_k^f - \bar{p}_f t\right)\right|  =  |I_k^f - \bar{p}_f|
    & \leq \max\{|0 - \bar{p}_f|, |1 - \bar{p}_f|\} = \bar{p}_f,
\end{align*}
for any $t \geq 1$, by the Azuma-Hoeffding inequality, we have for any 
positive constant $c$, $\P\left[ \sum^{t-1}_{k=0} I_k^f - p_f t \leq -c \right] \leq \exp\left( -\tfrac{c^2}{2 p_f^2 t} \right)$. 
Setting $c = (p_f - \bar{p}_f)T$ and subtracting 1 from both sides yields the second result.

Similarly, 
at the negative curvature steps, for all $k \in \N{}$ 
\begin{align*} 
        \P_{\hat{\xi}_k^{(0)}}\left[ \hat{E}_k > s | \hat{\mathcal{F}}_{k-1}' \right] \leq \exp (a(\epsilon_f - s)), \quad \P_{\hat{\xi}_k^{(0\pm)}} \left[ \hat{E}_k^\pm > s | \hat{\mathcal{F}}_{k-1}' \right] \leq \exp (a(\epsilon_f - s)) \nonumber
\end{align*}
Consequently, for all $k \in \N{}$,
\begin{align} \label{eq.prob_bound_max_hatEk}
        \P_{ \hat{\xi}_k^{(0\pm)}} \left[ \max\{\hat{E}_k^+, \hat{E}_k^-\} > s | \hat{\mathcal{F}}_{k-1}' \right] & \leq \P_{\hat{\xi}_k^{(0\pm)}} \left[ \hat{E}_k^+ > s \text{ or } \hat{E}_k^- > s | \hat{\mathcal{F}}_{k-1}' \right]\nonumber\\
        & \leq \P_{\hat{\xi}_k^{(0\pm)}} \left[ \hat{E}_k^+ > s | \hat{\mathcal{F}}_{k-1}' \right] + \P_{\hat{\xi}_k^{(0\pm)}} \left[ \hat{E}_k^- > s | \hat{\mathcal{F}}_{k-1}' \right] \nonumber\\
        & = \P_{\hat{\xi}_k^{(0+)}} \left[ \hat{E}_k^+ > s | \hat{\mathcal{F}}_{k-1}' \right] + \P_{\hat{\xi}_k^{(0-)}} \left[ \hat{E}_k^- > s | \hat{\mathcal{F}}_{k-1}' \right] \nonumber\\
        & \leq 2 \exp(a(\epsilon_f - s)) = \exp \left(a\left(\epsilon_f +\tfrac{1}{a} \log 2 - s \right) \right),
\end{align}
and as a result, by the definition of $I_k^f$, Oracle \ref{def.0th-order_orc}.\ref{orc.0th-order_subexp_orc}, and $p_f$, it follows that
\begin{align*}
    \P \left[ \hat{I}_k^f = 0| \mathcal{F}_{k-1}' \right] & = \P \left[\hat{E}_k + \max\{\hat{E}_k^+, \hat{E}_k^-\} > e_f| \hat{\mathcal{F}}_{k-1}' \right] \leq 3 \exp (a(\epsilon_f - \tfrac{e_f}{2})) = 1 - p_f.
\end{align*}
\eproof

The analysis of subexponential random variables draws on \cite[Chapter 2]{vershynin2018high}; however, those results are insufficient for our purposes, as thea distribution considered here is parameterized by two quantities, $a$ and $\epsilon_f$, rather than one. We therefore invoke a stronger result from \cite{cao2023first}, stated in the following proposition.
\bproposition \label{prop.subexp}
Let $X$ be a random variable such that for some $a > 0$ and $b \geq 0$, $\P[|X| \geq s] \leq \exp(a(b - s))$, for all $\ s > 0$.
Then, it follows that $\E \{\exp(\lambda|X|)\} \leq \tfrac{1}{1 - \lambda/a} \exp(\lambda b)$ for all $\lambda \in [0, a)$.
\eproposition

Using Proposition \ref{prop.subexp}, we next establish a lemma that provides a high-probability upper bound on the total error caused by noisy function evaluations.
\blemma \label{lmm.prob.ef_lg}
For any $s \geq 0$ and $t \geq 1$, $\P\left[ \sum_{k=0}^{t-1} (E_k + E_k^+) \geq t(\tfrac{5}{a} + 2\epsilon_f) + s \right]  \leq \exp\left(-\tfrac{a}{4}s \right)$ and $\P\left[ \sum_{k=0}^{t-1} \left(\hat{E}_k + \max\{\hat{E}_k^+, \hat{E}_k^-\} \right) \geq t(\tfrac{5}{a} + 2\epsilon_f) + s \right] \leq \exp\left(-\tfrac{a}{4}s \right)$.
\elemma
\bproof
    For the descent step, by Proposition \ref{prop.subexp}, since \eqref{eq.prob_bound_Ek} holds for both function estimation errors $E_k$ and $E_k^+$, it follows that
    \begin{align*}
         \E\left[ \exp(2\lambda E_k) \right] &\leq \tfrac{1}{1-2\lambda/a} \exp(2\lambda \epsilon_f) 
        \quad \text{ and }\quad   \E\left[ \exp(2\lambda E_k^+) \right] \leq \tfrac{1}{1-2\lambda/a} \exp(2\lambda \epsilon_f), 
    \end{align*}
    for all $\lambda \in \left[0, \tfrac{a}{2}\right)$. 
    Then, by the Cauchy-Schwarz inequality, for all $\lambda \in \left[0, \tfrac{a}{2}\right)$, 
    \begin{equation} \label{eq.subexp.expectation}
        \begin{split}
            \E \left[ \exp(\lambda E_k + \lambda E_k^+ ) | \mathcal{F}_{k-1}' \right] & \leq \sqrt{\E \left[ \exp(2\lambda E_k) | \mathcal{F}_{k-1}' \right] \E \left[ \exp(2\lambda E_k^+) | \mathcal{F}_{k-1}' \right]} \\
            & \leq \tfrac{1}{1-2\lambda/a} \exp(2\lambda \epsilon_f).
        \end{split}
    \end{equation}
    By Markov’s inequality, for any $\lambda \in \left[0, \tfrac{a}{2}\right)$, $s \geq 0$, and positive integer $t$, we have
    \begin{equation} \label{eq.prob_bound_sum_Ek}
        \begin{split}
            \P\left[ \sum_{k=0}^{t-1} (E_k + E_k^+) \geq s \right] & = \P\left[ \exp\left( \lambda \sum_{k=0}^{t-1} (E_k + E_k^+) \right) \geq \exp(\lambda s) \right] \\
            & \leq e^{-\lambda s}  \E \left[ \exp\left( \lambda \sum_{k=0}^{t-1} (E_k + E_k^+) \right) \right] \\
            & \leq e^{-\lambda s} \left( \tfrac{1}{1-2\lambda/a} \exp(2\lambda \epsilon_f) \right)^t,
        \end{split}
    \end{equation}
    where the last inequality can be proven by induction as follows. First, this inequality holds for $t=1$ due to \eqref{eq.subexp.expectation}. Now if it holds for any positive integer $t$, then for $t+1$,
    \begin{align*}
        & \quad\, e^{-\lambda s}  \E \left[ \exp\left( \lambda \sum_{k=0}^{(t+1)-1} (E_k + E_k^+) \right) \right] \\ 
        & = e^{-\lambda s}  \E \left[ \exp\left( \lambda \sum_{k=0}^{t-1} (E_k + E_k^+) \right)  \E_{\xi_t^{(0)}, \xi_t^{(0+)}} \left[ \exp\left( \lambda (E_t + E_t^+) \right) | \mathcal{F}_{t-1}' \right] \right] \\
        & \leq e^{-\lambda s}  \E \left[ \exp\left( \lambda \sum_{k=0}^{t-1} (E_k + E_k^+) \right) \tfrac{1}{1-2\lambda/a}\exp(2\lambda \epsilon_f) \right]  \leq e^{-\lambda s} \left( \tfrac{1}{1-2\lambda/a} \exp(2\lambda \epsilon_f) \right)^t
    \end{align*}
    where the first inequality holds by \eqref{eq.subexp.expectation} and the second inequality holds by the induction hypothesis. This completes the induction proof. 
    For ease of exposition, we use $1/(1 - x) \leq \exp(2x)$ for all $x \in [0, 1/2]$ to simplify the result
    \begin{align*}
        \P\left[ \sum_{k=0}^{t-1} (E_k + E_k^+) \geq s \right] \leq e^{-\lambda s} \left[ \exp(\tfrac{4\lambda}{a}) \exp(2\lambda \epsilon_f) \right]^t \leq \exp(\lambda[t (\tfrac{5}{a} + 2 \epsilon_f ) - s]),
    \end{align*}
    where $\lambda \in \left[0, \tfrac{a}{4}\right]$. The right-hand side is only less than or equal to 1 when $s \geq t (\tfrac{5}{a} + 2 \epsilon_f)$. Replacing $s$ with $t (\tfrac{5}{a} + 2 \epsilon_f) + s$ and choosing $\lambda = \tfrac{a}{4}$, yields the first result. 

    For the negative curvature step, by Proposition \ref{prop.subexp} and inequalities \eqref{eq.prob_bound_Ek}-\eqref{eq.prob_bound_max_hatEk}, 
    \begin{align*}
        \E\left[ \exp(2\lambda \hat{E}_k) \right] &\leq \tfrac{1}{1-2\lambda/a} \exp(2\lambda \epsilon_f),\\
        \text{and } \quad \E\left[ \exp(2\lambda \max\{\hat{E}_k^+, \hat{E}_k^-\}) \right] &\leq \tfrac{1}{1-2\lambda/a} \exp(2\lambda (\epsilon_f + \tfrac{1}{a}\log 2)),
    \end{align*}
    for all $\lambda \in \left[0, \tfrac{a}{2}\right)$. Following the same procedure as \eqref{eq.prob_bound_sum_Ek}, 
    \begin{align*}
        \P\left[ \sum_{k=0}^{t-1} \left(\hat{E}_k + \max\{\hat{E}_k^+, \hat{E}_k^-\} \right) \geq s \right] &\leq e^{-\lambda s} \left( \tfrac{1}{1-2\lambda/a} \exp\left(\lambda(2 \epsilon_f+ \tfrac{1}{a} \log 2) \right) \right)^t \\
        &\leq 
        \exp(\lambda[t (\tfrac{5}{a} + 2 \epsilon_f ) - s]).
    \end{align*}
    The right-hand side is only less than or equal to 1 when $s \geq t (\tfrac{5}{a} + 2 \epsilon_f)$. Replacing $s$ with $t (\tfrac{5}{a} + 2 \epsilon_f) + s$ and choosing $\lambda = \tfrac{a}{4}$, yields the second result. 
\eproof

We introduce the model accuracy assumptions (conditions on the oracle success probabilities) required to derive high-probability iteration-complexity bounds.

\bassumption \label{asm.prob}
We make the following assumptions with regards to the probabilities ($p_f$, $p_g$, and $p_H$) in the probabilistic oracles:
\begin{itemize}[leftmargin=0.5cm]
    \item[$(i)$] $\P \left[ I_k^f I_k^g = 1 | \mathcal{F}_{k-1} \right] \geq p_f p_g \ $ and $\ \P \left[ \hat{I}_k^f I_k^H = 1 | \mathcal{F}_{k-1}, \hat{x}_k \right] \geq p_f p_H$,
    \item[$(ii)$] $\P \left[I_k^f I_k^g\hat{I}_k^f I_k^H = 1 | \mathcal{F}_{k-1} \right] \geq p_f^2 p_g p_H$, $\quad \text{and} \qquad (iii) \ $ $p_f^2 p_g p_H + p_fp_g + p_fp_H-2 > 0$.
\end{itemize}
\eassumption
\bremark
The bounds in $(i)$-$(ii)$ imply that the processes $\sum_{k=0}^{t-1} I_k^f I_k^g - p_f p_g t$, $\sum_{k=0}^{t-1} \hat{I}_k^f I_k^H - p_f p_H t$, and $\sum_{k=0}^{t-1} I_k^f I_k^g \hat{I}_k^f I_k^H - p_f^2 p_g p_H t$ are submartingales. When $p_f = 1$, we have $I_k^f = \hat{I}_k^f = 1$ almost surely and Assumption~\ref{asm.prob} reduces to the corresponding condition used in the bounded-noise setting. Moreover, if the algorithm performs only descent steps, 
then $\hat{I}_k^{f}$ and $I_k^{H}$ can be omitted, and Assumption~\ref{asm.prob} simplifies to $\mathbb{P}\left[I_k^{f} I_k^{g}=1 | \mathcal{F}_{k-1}\right] \geq p_f p_g \eqqcolon p_{fg} > \tfrac12$, 
which matches the assumption in \cite{jin2021high}. Related combined-probability conditions of the form in $(ii)$ also appear in \cite{cao2023first}.
\eremark

The following lemma follows from Assumption~\ref{asm.prob}(ii) and the Azuma–Hoeffding inequality \cite{azuma1967weighted}. Its proof is similar to that of Lemma~\ref{lmm.f.accuracy} and is omitted for brevity.

\blemma \label{lmm.azuma}
For all positive integers $t$, and any $\bar{p}_f \in [0, p_f)$,  $\bar{p}_g \in [0, p_g)$, and $\bar{p}_H \in [0, p_H)$, $\P \left[ \sum_{k=0}^{t-1} I_k^f I_k^g > \bar{p}_f \bar{p}_g t\right] \geq 1- \exp \left( - \tfrac{(p_f p_g - \bar{p}_f \bar{p}_g)^2}{2 p_f^2 p_g^2} t \right)$, \\$\P \left[ \sum_{k=0}^{t-1} \hat{I}_k^f I_k^H > \bar{p}_f \bar{p}_H t\right] \geq 1- \exp \left( - \tfrac{(p_f p_H - \bar{p}_f \bar{p}_H)^2}{2 p_f^2 p_H^2} t \right)$, and \\$\P \left[ \sum_{k=0}^{t-1} I_k^f I_k^g \hat{I}_k^f I_k^H > \bar{p}_f^2 \bar{p}_g \bar{p}_H t\right] \geq 1- \exp \left( - \tfrac{(p_f^2 p_g p_H - \bar{p}_f^2 \bar{p}_g \bar{p}_H)^2}{2 p_f^4 p_g^2 p_H^2} t \right)$.
\elemma
The results of Lemmas \ref{lmm.f.accuracy}, \ref{lmm.prob.ef_lg}, and \ref{lmm.azuma}, serve as central building blocks for the subexponential noise analysis in Section \ref{subsec.analysis.subexp}.

\section{Convergence Analysis} \label{sec.analysis}

In this section, we present a comprehensive convergence and complexity analysis for Algorithm~\ref{alg.theory}. For convenience in the analysis, we introduce the following indicator random variables:
\begin{align*}
    \Omega_k^g & \coloneqq \mathbbm{1}\left\{ \|g_k\|_2 \geq c_g \bar{\epsilon}_g \right\}, \; \Omega_k^H \coloneqq \mathbbm{1}\left\{ \lambda_k < - c_H \max\{\bar{\epsilon}_H, \bar{\epsilon}_{\lambda} \right\}\},\\
    \Theta_k^g & \coloneqq \mathbbm{1}\{\text{descent step is successful, i.e., } \alpha_{k+1} = \tau^{-1} \alpha_k\}, \\
    \text{and } \quad \Theta_k^H & \coloneqq \mathbbm{1}\{\text{negative curvature step is successful, i.e., } \beta_{k+1} = \tau^{-1} \beta_k\}.
\end{align*}
To formalize the iteration complexity, we define the stopping time
\begin{equation} \label{def.stopping_time}
    N_{\bar{\epsilon}_g, \bar{\epsilon}_H, \bar{\epsilon}_{\lambda}} \coloneqq \min \left\{k : \|\nabla f(x_k)\|_2 \leq \bar{\epsilon}_g, \, \lambda_{\min}(\nabla^2 f(x_k)) \geq - \max\{\bar{\epsilon}_{\lambda}, \bar{\epsilon}_{H}\} \right\},
\end{equation}
i.e., the number of iterations required to reach an ($\bar{\epsilon}_g,\bar{\epsilon}_{\lambda},\bar{\epsilon}_H$)-stationary point. 

\subsection{Convergence analysis: Bounded noise case}
In this subsection, we analyze convergence under Oracle~\ref{def.0th-order_orc}.\ref{orc.0th-order_bounded}, the bounded function evaluation oracle. Specifically, the realizations of the function estimate errors satisfy
\begin{equation} \label{eq.bounded_anal.bounded_noise}
    e_k \leq \epsilon_f, \ e_k^+ \leq \epsilon_f, \ \hat{e}_k \leq \epsilon_f, \ \hat{e}_k^{\pm} \leq \epsilon_f.
\end{equation}
The line search $c_d$, $c_p$ and early termination $c_g$ and $c_H$ parameters are set to satisfy
\begin{equation*}
    \begin{split}
        & 0 < c_d < \tfrac{1}{2} + \tfrac{(1-\eta) (1-\kappa_g)}{2(1+\kappa_g)}, \quad  0 < c_g < 1-\kappa_g  \leq 1, \quad 0 < \eta < 1, \\
        \text{and} \qquad & 0 < c_p < \tfrac{\gamma-\kappa_H}{2\gamma} = \tfrac{1}{2}- \tfrac{\kappa_H}{2\gamma}, \quad 0 \leq \kappa_H < \gamma \leq 1, \quad 0 < c_H < 1-\kappa_{\lambda} \leq 1.
    \end{split}
\end{equation*}

These bounds quantify how oracle accuracy influences the admissible parameter ranges. For the gradient step, as $\kappa_g \to 1$, the upper bounds on $c_d$ and $c_g$ decrease, necessitating more conservative acceptance and termination thresholds under less accurate gradient estimates. For the negative curvature step, increasing $\kappa_{\lambda}$ drives the allowable range of $c_H$ toward zero, imposing stricter requirements on curvature-based termination, while $\kappa_H \to 1$ forces $\gamma \to 1$, enforcing a stronger criterion for accepting a negative curvature direction (extreme case $\gamma = 1$, direction must coincide with eigenvector). In the noiseless setting, i.e., $\epsilon_g = \kappa_g = 0$ and $\epsilon_H = \kappa_H = \epsilon_{\lambda} = \kappa_{\lambda} = 0$, the requirements reduce to $0 < c_d < 1-\tfrac{\eta}{2}$, $0 < c_g, c_H < 1$, $0 < c_p < \tfrac{1}{2}$, and $0 < \gamma \leq 1$, which coincide with the bounds in the deterministic setting.

We next establish properties of the stochastic process generated by Algorithm~\ref{alg.theory} when the function, gradient, and Hessian estimates are computed via Oracles~\ref{def.0th-order_orc}.\ref{orc.0th-order_bounded}, \ref{def.1st-order_orc}, and \ref{def.2nd-order_orc}, respectively.
\blemma \label{lmm.bounded.stocproc}
    Let $e_f \geq 2\epsilon_f$, and 
    \begin{align} \label{eq.value.eps}
        \bar{\epsilon}_g \geq \tfrac{2\epsilon_g}{\min\{\eta c_g (1 - \kappa_g), 1-\kappa_g-c_g\}}, 
        \ \bar{\epsilon}_H \geq \tfrac{\epsilon_{H}}{c_H} \sqrt{\tfrac{2}{\delta(\gamma- \kappa_H -2 c_p \gamma)}}, \ \text{ and } \ \bar{\epsilon}_{\lambda} \geq \tfrac{\epsilon_{\lambda}}{1-\kappa_{\lambda}-c_H}.
    \end{align}
    Then, there exist constants $\bar{\alpha}, \bar{\beta} > 0$ 
    \begin{equation} \label{def.bar_alpha_beta}
        \bar{\alpha} \coloneqq \tfrac{2 \left( 1/(1+\kappa_g') -c_d \right)}{L_g}, \quad \text{and} \quad \bar{\beta} \coloneqq \tfrac{3(\gamma-\kappa_{H}- 2c_p \gamma)}{2\delta L_H},
    \end{equation}
    where $0 \leq \kappa_g' < 1$, 
    and nondecreasing functions $h_d(\cdot), h_p(\cdot): \mathbb{R} \rightarrow \mathbb{R}$, satisfying $h_d(\alpha), h_p(\beta) > 0$ for any $\alpha, \beta > 0$, such that for any realization of Algorithm \ref{alg.theory} the following properties hold for all $k < N_{\bar{\epsilon}_g, \bar{\epsilon}_H, \bar{\epsilon}_{\lambda}}$:
    \begin{enumerate}[leftmargin=0.5cm]
        \item [$(i)$] If the gradient estimate at the descent step of iteration $k$ is accurate (i.e., $I_k^g = 1$) and $\|g_k\|_2 \geq c_g \bar{\epsilon}_g$ (i.e., $\Omega_k^g = 1$), then for $\alpha_k \leq \bar{\alpha}$, the iteration (descent step) is successful ($\Theta_k^g = 1$), which implies $\alpha_{k+1} = \tau^{-1} \alpha_k$;
        \item [$(ii)$] If the Hessian estimate at the negative curvature step of iteration $k$ is accurate (i.e., $I_k^H = 1$) and $\lambda_k \leq -c_H \max \{\bar{\epsilon}_H, \bar{\epsilon}_{\lambda}\}$ (i.e., $\Omega_k^H = 1$), then for $\beta_k \leq \bar{\beta}$, the iteration (negative curvature) is successful ($\Theta_k^H = 1$), which implies $\beta_{k+1} = \tau^{-1} \beta_k$;
        \item [$(iii)$] When $I_k^g \Omega_k^g \Theta_k^g = 1$, then $f(x_{k+1}) \leq f(x_k) - h_d(\alpha_k) + 4e_f$;
        \item [$(iv)$] When $I_k^H \Omega_k^H \Theta_k^H = 1$, then $f(x_{k+1}) \leq f(x_k) - h_p(\beta_k) + 4e_f$;
        \item [$(v)$] $I_k^g(1-\Omega_k^g) I_k^H(1-\Omega_k^H) = 0$.
    \end{enumerate}
    In summary,
    \begin{align} \label{eq.summary_bounded}
        f(x_{k+1}) \leq \begin{cases}
            f(x_k) - h_d(\alpha_k) + 4e_f, & \text{if } I_k^g \Omega_k^g \Theta_k^g = 1, \\
            f(x_k) - h_p(\beta_k) + 4e_f, & \text{if } I_k^H \Omega_k^H \Theta_k^H = 1, \\
            f(x_k) + 4e_f, & \text{otherwise}.
        \end{cases}
    \end{align}
\elemma 
\bproof 
    We start with the descent step and consider the case where $I_k^g \Omega_k^g = 1$ (gradient is accurate and $\|g_k\|_2 > c_g \bar{\epsilon}_g$). By Oracle~\ref{def.1st-order_orc}, \eqref{eq.value.eps}, and the termination condition in Algorithm~\ref{alg.theory}, we have $\|\nabla f(x_k) - g_k\|_2 \leq \epsilon_g + \kappa_g \|\nabla f(x_k) \|_2$ and $\|g_k\|_2 \geq c_g \bar{\epsilon}_g \geq \tfrac{2c_g \epsilon_g}{\eta c_g(1-\kappa_g)} = \tfrac{2\epsilon_g}{\eta (1-\kappa_g)}$, and as a result 
    $\|\nabla f(x_k) - g_k\|_2 \leq  \kappa_g' \|\nabla f(x_k)\|_2$, where $\kappa_g' \coloneqq 1- \tfrac{2(1-\eta) (1-\kappa_g)}{1+\kappa_g + (1-\eta) (1-\kappa_g)} < 1$, from which it follows that 
    \begin{align}
        -\nabla f(x_k)^{\mathsf{T}}g_k \leq - \tfrac{1}{1 + \kappa_g'} \|g_k\|_2^2, \quad -\nabla f(x_k)^{\mathsf{T}} g_k \leq -(1 - \kappa_g') \|\nabla f(x_k)\|_2^2. \label{eq.bounded.g_relative_bound}
    \end{align}
    Then, by Assumption~\ref{asm.f.general}, $d_k = - g_k$, and \eqref{eq.bounded.g_relative_bound},
    \begin{align*}
            f(x_k + \alpha_kd_k) & \leq f(x_k) + \alpha_k \nabla f(x_k)^{\mathsf{T}} d_k + \tfrac{L_g}{2} \alpha_k^2 \|d_k\|_2^2 \\
            & \leq f(x_k) - \tfrac{\alpha_k}{1 + \kappa_g'} \|g_k\|_2^2 + \tfrac{L_g}{2} \alpha_k^2 \|g_k\|_2^2.
    \end{align*}
    Thus, as long as $\alpha_k \leq \bar{\alpha} = \tfrac{2 \left( 1/(1+\kappa_g') -c_d \right)}{L_g}$
    the inequality can be further bounded by
    \begin{equation} \label{eq.bounded.fk+bound}
        f(x_k + \alpha_kd_k) \leq f(x_k) - c_d \alpha_k \|g_k\|_2^2 = f(x_k) + c_d \alpha_k g_k^{\mathsf{T}} d_k.
    \end{equation}
    Invoking the bounded noise condition \eqref{eq.bounded_anal.bounded_noise} and $e_f \geq 2\epsilon_f$, 
    \begin{equation} \label{eq.bounded.Fk+bound}
        \begin{split}
            F(x_k + \alpha_kd_k, \xi_k^{(0+)})  \leq f(x_k + \alpha_kd_k) + e_k^+ 
            & \leq f(x_k) + c_d \alpha_k g_k^{\mathsf{T}} d_k + \epsilon_f \\
            & \leq F(x_k, \xi_k^{(0)}) + c_d \alpha_k g_k^{\mathsf{T}} d_k + e_f,
        \end{split}
    \end{equation}
    which implies that the descent step at iterate $k$ is successful ($\Theta_k^g = 1$) provided the step size is sufficiently small (part $(i)$ of Lemma~\ref{lmm.bounded.stocproc}).
    
    When the step is also successful, that is $I_k^g \Omega_k^g \Theta_k^g = 1$, according to the procedure in Algorithm \ref{alg.theory} and the noise bound \eqref{eq.bounded_anal.bounded_noise}, it follows that
    \begin{align*}
            F(x_k+ \alpha_kd_k, \xi_k^{(0+)})  \leq 
            F(x_k, \xi_k^{(0)}) - c_d \alpha_k \|g_k\|_2^2 + e_f  &\leq F(x_k, \xi_k^{(0)}) - c_d \alpha_k c_g^2 \bar{\epsilon}_g^2 + e_f, 
    \end{align*}
    where $h_d(\alpha) \coloneqq c_g^2 \bar{\epsilon}_g^2 \alpha$. Again, invoking the bounded noise condition  and $e_f \geq 2\epsilon_f$,
    \begin{equation} \label{eq.bounded.dc_bound_good}
        \begin{split}
            f(\hat{x}_k) = f(x_k+ \alpha_kd_k) & \leq F(x_k+ \alpha_k d_k, \xi_k^{(0+)}) + e_k^+ \\
            & \leq F(x_k, \xi_k^{(0)}) - h_d(\alpha_k) + e_f + \epsilon_f \\
            & \leq f(x_k) - h_d(\alpha_k) + 2e_f,
        \end{split}
    \end{equation}
    establishing a decrease in the true function.
    
    If instead, $I_k^g (1-\Omega_k^g) = 1$, then $\|\nabla f(x_k)\|_2 \leq \|g_k\|_2 + \|\nabla f(x_k)-g_k \|_2 < c_g \bar{\epsilon}_g + \epsilon_g + \kappa_g \|\nabla f(x_k)\|_2$ 
    from which it follows that
    \begin{equation} \label{eq.bounded.stop_cond_k}
        \|\nabla f(x_k)\|_2  \leq \tfrac{c_g \bar{\epsilon}_g + \epsilon_g}{1-\kappa_g} \leq \tfrac{2c_g + 1-\kappa_g-c_g}{2(1-\kappa_g)} \bar{\epsilon}_g < \bar{\epsilon}_g,
    \end{equation}
    hence the first part of the stopping time definition \eqref{def.stopping_time} holds.

    When $\Theta_k^g = 0$ (the descent step is unsuccessful), the trial point is rejected and $\hat{x}_k = x_k$, $f(\hat{x}_k) = f(x_k)$. Otherwise, $\Theta_k^g = 1$ ensures that the potential increase in the function value is always bounded
    \begin{equation} \label{eq.bounded.dc_bound_bad}
        \begin{split}
            f(x_k+ \alpha_kd_k) & \leq F(x_k+ \alpha_kd_k, \xi_k^{(0+)}) + e_k^+ \\
            & \leq F(x_k, \xi_k^{(0)}) + c_d \alpha_k g_k^{\mathsf{T}} d_k + e_f + \epsilon_f \\
            & \leq F(x_k, \xi_k^{(0)}) + e_f + \epsilon_f  
            \leq f(x_k) + 2e_f
        \end{split}
    \end{equation}
    We should note that this bound also applies to the case where $\Theta_k^g = 0$.

    Next, at $\hat{x}_k \leftarrow x_k + \alpha_k d_k$, we consider the negative curvature step. When $I_k^H \Omega_k^H = 1$ (the Hessian is accurate and $\lambda_k < - c_H \max \{\bar{\epsilon}_H, \bar{\epsilon}_{\lambda}\}$), by \eqref{eq.value.eps} and \eqref{fom.q2} 
    \begin{equation} \label{eq.bounded.hessian_vector}
        \begin{split}
            \|(\nabla^2 f(\hat{x}_k) - H_k) q_k\|_2 & \leq \epsilon_H^2 + \kappa_H |\lambda_k| \|q_k\|_2 \\
            & \leq \tfrac{\delta (\gamma-\kappa_{H} - 2 c_p \gamma) (c_H \bar{\epsilon}_H)^2}{2} + \kappa_H \delta |\lambda_k|^2 \\
            & < \tfrac{\delta (\gamma + \kappa_{H} - 2 c_p \gamma) |\lambda_k|^2}{2},
        \end{split}
    \end{equation}
    By the Lipschitz continuity of $\nabla^2 f$ (Assumption~\ref{asm.f.general}),
    \begin{equation*}
        f(\hat{x}_k \pm \beta_k q_k) \leq f(\hat{x}_k) \pm \beta_k \nabla f(\hat{x}_k)^{\mathsf{T}} q_k + \tfrac{\beta_k^2}{2} q_k^{\mathsf{T}} \nabla^2 f(\hat{x}_k) q_k + \tfrac{L_H}{6} \beta_k^3 \|q_k\|_2^3.
    \end{equation*}
    Since $\min \{\nabla f(\hat{x}_k)^{\mathsf{T}} q_k, -\nabla f(\hat{x}_k)^{\mathsf{T}} q_k\} \leq 0$, taking the minimum of the above
    \begin{equation} \label{eq.bounded.min_nc_upp}\
        \min \{f(\hat{x}_k +\beta_k q_k), f(\hat{x}_k-\beta_k q_k)\} \leq f(\hat{x}_k) + \tfrac{\beta_k^2}{2} q_k^{\mathsf{T}} \nabla^2 f(\hat{x}_k) q_k + \tfrac{L_H}{6} \beta_k^3 \|q_k\|_2^3.
    \end{equation}
    Notice that $q_k^{\mathsf{T}} \nabla^2 f(\hat{x}_k) q_k = q_k^{\mathsf{T}} (\nabla^2 f(\hat{x}_k)-H_k) q_k + q_k^{\mathsf{T}} H_k q_k \leq \|q_k\|_2 \|(\nabla^2 f(\hat{x}_k)-H_k) q_k\|_2 + q_k^{\mathsf{T}} H_k q_k$. By \eqref{eq.bounded.hessian_vector} for any $\beta_k \leq \bar{\beta}$, it follows that
    \begin{align*}
        \tfrac{\beta_k^2}{2} \|q_k\|_2 \|(\nabla^2 f(\hat{x}_k)-H_k) q_k\|_2 + \tfrac{\beta_k^2}{2} q_k^{\mathsf{T}} H_k q_k + \tfrac{L_H}{6} \beta_k^3 \|q_k\|_2^3 \leq c_p \beta_k^2 q_k^{\mathsf{T}} H_k q_k.
    \end{align*}
    Substituting the above bound into \eqref{eq.bounded.min_nc_upp}, when $I_k^g \Omega_k^g = 1$ and $\beta_k \leq \bar{\beta}$, 
    \begin{align*}
        \min \left\{f(\hat{x}_k +\beta_k q_k), f(\hat{x}_k-\beta_k q_k) \right\} \leq f(\hat{x}_k) + c_p \beta_k^2 q_k^{\mathsf{T}} H_k q_k,
    \end{align*}
    and consequently, by \eqref{eq.bounded_anal.bounded_noise} and $e_f \geq 2\epsilon_f$, it follows that
    \begin{align*}
        & \quad \, \min \left\{F(\hat{x}_k +\beta_k q_k, \hat{\xi}_k^{(0+)}), F(\hat{x}_k-\beta_k q_k, \hat{\xi}_k^{(0-)})\right\} \\
        & \leq \min \left\{f(\hat{x}_k +\beta_k q_k) + \hat{e}_k^+, f(\hat{x}_k-\beta_k q_k) + \hat{e}_k^- \right\} \\
        & \leq \min \left\{f(\hat{x}_k +\beta_k q_k), f(\hat{x}_k-\beta_k q_k)\right\} + \max\{\hat{e}_k^+, \hat{e}_k^-\} \\
        & \leq f(\hat{x}_k) + c_p \beta_k^2 q_k^{\mathsf{T}} H_k q_k + \epsilon_f \leq F(\hat{x}_k, \hat{\xi}_k^{(0)}) + c_p \beta_k^2 q_k^{\mathsf{T}} H_k q_k + e_f.
    \end{align*}
    Therefore, $\Theta_k^H = 1$ provided the step size is sufficiently small (part $(ii)$ of Lemma~\ref{lmm.bounded.stocproc}). 
    
    When $I_k^H \Omega_k^H \Theta_k^H = 1$, the step is successful and $\lambda_k < - c_H \max \{\bar{\epsilon}_H, \bar{\epsilon}_{\lambda}\}$, and 
    \begin{align} 
            & \quad \,  \min \left\{F(\hat{x}_k + \beta_k q_k, \hat{\xi}_k^{(0+)}), F(\hat{x}_k - \beta_k q_k, \hat{\xi}_k^{(0-)})\right\} \notag \\
            & \leq F(\hat{x}_k, \hat{\xi}_k^{(0)}) + c_p \beta_k^2 p_k^{\mathsf{T}} H_k p_k + e_f \notag\\
            & \leq F(\hat{x}_k, \hat{\xi}_k^{(0)}) - c_p \gamma \delta^2 \beta_k^2 |\lambda_k|^3 + e_f \label{eq.bounded.nc_bound_decrease}\\ 
            & \leq F(\hat{x}_k, \hat{\xi}_k^{(0)}) - \beta_k^2 c_p \gamma \delta^2 c_H^3 (\max \{\bar{\epsilon}_H, \bar{\epsilon}_{\lambda}\})^3 + e_f, \notag  
    \end{align}
    where $h_p(\beta) \coloneqq c_p \gamma \delta^2 c_H^3 (\max \{\bar{\epsilon}_H, \bar{\epsilon}_{\lambda}\})^3 \beta^2$. 
    If $F(\hat{x}_k \pm \beta_k q_k, \hat{\xi}_k^{(0+)}) = \min\{F(\hat{x}_k + \beta_k q_k, \hat{\xi}_k^{(0+)}), F(\hat{x}_k - \beta_k q_k, \hat{\xi}_k^{(0-)})\}$, by the definition of $p_k$ and \eqref{eq.bounded.nc_bound_decrease}
    \begin{equation} \label{eq.bounded.nc_bound_+q}
        \begin{split}
            f(\hat{x}_k + \beta_k p_k) &= f(\hat{x}_k \pm \beta_k q_k) \\
            & \leq F(\hat{x}_k \pm \beta_k q_k, \hat{\xi}_k^{(0+)}) + \hat{e}_k^+ \\
            & = \min \{F(\hat{x}_k + \beta_k q_k, \hat{\xi}_k^{(0+)}), F(\hat{x}_k - \beta_k q_k, \hat{\xi}_k^{(0-)})\} + \hat{e}_k^+ \\
            & \leq F(\hat{x}_k, \hat{\xi}_k^{(0)}) - h_p(\beta_k) + e_f + \max\{\hat{e}_k^+, \hat{e}_k^-\}.
        \end{split}
    \end{equation}
    Combining the two cases 
    \begin{equation}\label{eq.bounded.nc_bound_good}
        \begin{split}
            f(\hat{x}_k + \beta_k p_k) & \leq F(\hat{x}_k, \hat{\xi}_k^{(0)}) - h_p(\beta_k) + e_f + \max\{\hat{e}_k^+, \hat{e}_k^-\} \\
            & \leq f(\hat{x}_k) + \epsilon_f - h_p(\beta_k) + e_f + \epsilon_f 
            \leq f(\hat{x}_k) - h_p(\beta_k) + 2e_f. 
         \end{split}
    \end{equation}
    
    When $I_k^H (1-\Omega_k^H) = 1$, we have $\lambda_k \geq -c_H \max \{\bar{\epsilon}_H, \bar{\epsilon}_{\lambda}\}$ and by \eqref{eq.value.eps}
    \begin{equation*}
        \begin{split}
            |\lambda_{\min} (\nabla^2 f(\hat{x}_k)) - \lambda_k| & \leq \epsilon_{\lambda} + \kappa_{\lambda} |\lambda_{\min} (\nabla^2 f(\hat{x}_k))| \\
            & \leq (1-\kappa_{\lambda}-c_H) \bar{\epsilon}_{\lambda}
            + \kappa_{\lambda} |\lambda_{\min} (\nabla^2 f(\hat{x}_k))| \\
            & \leq (1-\kappa_{\lambda}-c_H) \max \{\bar{\epsilon}_H, \bar{\epsilon}_{\lambda}\} + \kappa_{\lambda} |\lambda_{\min} (\nabla^2 f(\hat{x}_k))|,
        \end{split}
    \end{equation*}
    from which it follows that
    \begin{equation} \label{eq.bounded.stop_cond_hatk}
        \begin{split}
            |\lambda_{\min} (\nabla^2 f(\hat{x}_k)) - \lambda_k| & \geq \lambda_k - \lambda_{\min} (\nabla^2 f(\hat{x}_k)) \\
            & \geq -c_H \max \{\bar{\epsilon}_H, \bar{\epsilon}_{\lambda}\} - \lambda_{\min} (\nabla^2 f(\hat{x}_k)),
        \end{split}
    \end{equation}
    and $\lambda_{\min} (\nabla^2 f(\hat{x}_k)) \geq -\max \{\bar{\epsilon}_H, \bar{\epsilon}_{\lambda}\}$.
    
    When $\Theta_k^H = 0$ (the negative curvature step  is unsuccessful), the trial point is rejected and $x_{k+1} = \hat{x}_k$ and $f(x_{k+1}) = f(\hat{x}_k)$. Otherwise, $\Theta_k^H = 1$ ensures that
    \begin{align} 
            f(x_{k+1}) = f(\hat{x}_k+ \beta_kp_k) & \leq F(\hat{x}_k+ \beta_kp_k, \hat{\xi}_k^{(0\pm)}) + \max\{\hat{e}_k^+, \hat{e}_k^-\} \notag \\
            & \leq F(\hat{x}_k, \hat{\xi}_k^{(0)}) + c_p \beta_k^2 p_k^{\mathsf{T}} H_k p_k + e_f + \max\{\hat{e}_k^+, \hat{e}_k^-\} \notag \\
            & \leq F(\hat{x}_k, \hat{\xi}_k^{(0)}) - \gamma c_p \beta_k^2 |\lambda_k| \|p_k\|_2^2 + e_f + \max\{\hat{e}_k^+, \hat{e}_k^-\} \label{eq.bounded.nc_bound_bad}\\
            & \leq F(\hat{x}_k, \hat{\xi}_k^{(0)})  + e_f + \max\{\hat{e}_k^+, \hat{e}_k^-\} \notag\\
            & \leq f(\hat{x}_k) + e_f + 2\epsilon_f  \leq f(\hat{x}_k) + 2e_f. \notag
    \end{align}
    In both cases, $\Theta_k^H = 0$ and $\Theta_k^H = 1$, we have $f(x_{k+1}) \leq f(\hat{x}_k) + 2e_f$.
    
    Combining \eqref{eq.bounded.dc_bound_good}, \eqref{eq.bounded.nc_bound_bad} and $I_k^g \Omega_k^g \Theta_k^g = 1$ yields the first result in \eqref{eq.summary_bounded} (part $(iii)$ of Lemma~\ref{lmm.bounded.stocproc}). Combining \eqref{eq.bounded.dc_bound_bad}, \eqref{eq.bounded.nc_bound_good} and $I_k^H \Omega_k^H \Theta_k^H = 1$ yields the second result in \eqref{eq.summary_bounded} (part $(iv)$ of Lemma~\ref{lmm.bounded.stocproc}). The final case can be derived by \eqref{eq.bounded.dc_bound_bad} and \eqref{eq.bounded.nc_bound_bad}. 
    
    For $(v)$, by contradiction, if $I_k^g (1-\Omega_k^g) I_k^H (1-\Omega_k^H) = 1$, then $\hat{x}_k = x_k$ and using \eqref{eq.bounded.stop_cond_k} and \eqref{eq.bounded.stop_cond_hatk}, we can conclude that
    $$
        \|\nabla f(x_k)\|_2 < \bar{\epsilon}_g, \quad \lambda_{\min} (\nabla^2 f(x_k)) = \lambda_{\min} (\nabla^2 f(\hat{x}_k)) \geq -\max \{\bar{\epsilon}_H, \bar{\epsilon}_{\lambda}\},
    $$
    therefore $N_{\bar{\epsilon}_g, \bar{\epsilon}_H, \bar{\epsilon}_{\lambda}} \leq k$ which contradicts with the hypothesis that $k < N_{\bar{\epsilon}_g, \bar{\epsilon}_H, \bar{\epsilon}_{\lambda}}$.
\eproof

Lemma~\ref{lmm.bounded.stocproc} provides a fundamental bound on the change in the objective function value across iterations, accounting for the different cases that may arise in the algorithm. With the step size bounds $\bar{\alpha}$ and $\bar{\beta}$ specified in the lemma, we define the following random indicator variables to indicate if iteration $k$ has a large/small step.
\begin{align*}
    U_k^g = \begin{cases}
        1, & \min\{\alpha_k, \alpha_{k+1}\} \geq \bar{\alpha}, \Omega_k^g = 1 \\
        0, & \max\{\alpha_k, \alpha_{k+1}\} \leq \bar{\alpha}, \Omega_k^g = 1 \\
        0, & \text{if } \Omega_k^g = 0,
    \end{cases}, \ \  
    U_k^H = \begin{cases}
        1, & \min\{\beta_k, \beta_{k+1}\} \geq \bar{\beta}, \Omega_k^H = 1 \\
        0, & \max\{\beta_k, \beta_{k+1}\} \leq \bar{\beta}, \Omega_k^H = 1 \\
        0, & \Omega_k^H = 0.
    \end{cases}
\end{align*}
Without loss of generality, we assume that $\bar{\alpha} = \alpha_0 \tau^{m_{\alpha}}$ and $\bar{\beta} = \beta_0 \tau^{m_{\beta}}$ for some positive integers $m_{\alpha}$, $m_{\beta}$. In this case, it can be shown that for every step, $\alpha_k$ and $\beta_k$ is either a large step or a small step or a stable step (e.g., $\Omega_k^g$ or $\Omega_k^H = 0$). As an immediate deduction from properties $(i)$-$(ii)$ in Lemma \ref{lmm.bounded.stocproc}, it follows that 
\begin{equation} \label{eq.bounded.impossible_case}
    I_k^g \Omega_k^g (1-U_k^g) (1-\Theta_k^g) = 0, \quad \text{and} \quad I_k^H \Omega_k^H (1-U_k^H) (1-\Theta_k^H) = 0.
\end{equation}
The above implies that if the gradient or negative curvature information is sufficiently accurate and of sufficiently large magnitude, then any step size below certain threshold yields a successful step.
In Lemmas \ref{lmm.bounded.good_upper}-\ref{lmm.bounded.good_lower}, we provide bounds on the number of iterations for different cases.

\blemma \label{lmm.bounded.good_upper}
For any positive integer $t$, 
\begin{align*}
    \sum_{k=0}^{t-1} I_k^g \Omega_k^g \Theta_k^g U_k^g + \sum_{k=0}^{t-1} I_k^H \Omega_k^H \Theta_k^H U_k^H \leq \tfrac{f(x_0)-f(x_t)}{c_{\bar{\alpha}, \bar{\beta}}} + \tfrac{4e_f}{c_{\bar{\alpha}, \bar{\beta}}} t.
\end{align*}
where $c_{\bar{\alpha}, \bar{\beta}} \coloneqq \min\{h_d(\bar{\alpha}), h_p(\bar{\beta})\}$ and $\bar\alpha,\bar\beta$ are given in \eqref{def.bar_alpha_beta}.
\elemma
\bproof
Taking the step size indicator variables ($U_k^g$ and $U_k^H$) into consideration, one can derive a slight modification of \eqref{eq.summary_bounded} from Lemma \ref{lmm.bounded.stocproc}, i.e., 
\begin{equation*}
    f(x_{k+1}) \leq \begin{cases}
        f(x_k) - h_d(\bar{\alpha}) + 4e_f, & \text{if } I_k^g\Omega_k^g \Theta_k^g U_k^g = 1 \\
        f(x_k) - h_p(\bar{\beta}) + 4e_f, & \text{if } I_k^H\Omega_k^H \Theta_k^H U_k^H = 1 \\
        f(x_k) + 4e_f, & \text{otherwise,}
    \end{cases}
\end{equation*}
The above holds since $h_d(\cdot)$ and $h_p(\cdot)$ are non-decreasing functions on $\mathbb{R}_{\geq 0}$.
By summing up the inequalities from $k=0$ to $t-1$, it follows that
\begin{equation*}
    \begin{split}
        f(x_t) \leq f(x_0) -  h_d(\bar{\alpha})\sum_{k=0}^{t-1} I_k^g\Omega_k^g \Theta_k^g U_k^g - h_p(\bar{\beta})\sum_{k=0}^{t-1} I_k^H\Omega_k^H \Theta_k^H U_k^H   + 4t  e_f.
    \end{split}
\end{equation*}
Re-arranging the above and using the definition of $c_{\bar{\alpha}, \bar{\beta}}$ completes the proof.
\eproof

Lemma~\ref{lmm.bounded.good_upper} follows from Lemma~\ref{lmm.bounded.stocproc} and bounds the number of “good” iterations (accurate, successful, large step sizes). The next lemma bounds the number of iterations in which the step sizes $\alpha_k$ and $\beta_k$ are modified. When $\|g_k\|_2$ is sufficiently large ($\Omega_k^g = 1$), $\alpha_k$ increases on successful iterations ($\Theta_k^g = 1$) and decreases otherwise. Hence, the number of iterations with $U_k^g = 1$ in which $\alpha_k$ decreases is bounded by those in which it increases, plus the number required to reduce $\alpha_k$ from $\alpha_0$ to $\bar{\alpha}$. Similarly, the number of iterations with $\Omega_k^g(1-U_k^g)=1$ in which $\alpha_k$ increases is bounded by those in which it decreases. The following lemma formalizes these observations.
\blemma \label{lmm.bounded.bad_upper}
    For any positive integer $t$,
    \begin{align}
        \sum_{k=0}^{t-1} \Omega_k^g (1-U_k^g) \Theta_k^g + \sum_{k=0}^{t-1} \Omega_k^g U_k^g (1-\Theta_k^g) &\leq \tfrac{1}{2} \left(\sum_{k=0}^{t-1} \Omega_k^g + c_{\tau} \right), \label{eq.g.suc_bound_comb} \\
        \sum_{k=0}^{t-1} \Omega_k^H (1-U_k^H) \Theta_k^H + \sum_{k=0}^{t-1} \Omega_k^H U_k^H (1-\Theta_k^H) &\leq \tfrac{1}{2} \left(\sum_{k=0}^{t-1} \Omega_k^H + c_{\tau} \right), \label{eq.H.suc_bound_comb}
    \end{align}
    where $c_{\tau} \coloneqq \max\{\log_{\tau} \tfrac{\bar{\alpha}}{\alpha_0}, \log_{\tau} \tfrac{\bar{\beta}}{\beta_0},0\}$ and $\bar\alpha,\bar\beta$ are given in \eqref{def.bar_alpha_beta}.
\elemma
\bproof 
We first consider the bound for descent steps \eqref{eq.g.suc_bound_comb}. For all $\alpha_t$, $t \in \mathbb{N}$, the step size can be expressed as 
\begin{equation}\label{eq.stepgeneral}
    \begin{split}
        \alpha_t & = \alpha_0  \tau^{\sum_{k=0, \, \Omega_k^g = 1}^{t-1}  ((1-\Theta_k^g) -\Theta_k^g)} = \alpha_0  \tau^{\sum_{k=0}^{t-1} \Omega_k^g ((1-\Theta_k^g) -\Theta_k^g)}.
    \end{split}
\end{equation}
Let $s_{-1} \coloneqq 0$, $s_0 \coloneqq \min\{k \in \N{}: \alpha_k = \bar{\alpha}\}$, and $s_i \coloneqq \min\{k: k > s_{i-1}, \alpha_k = \bar{\alpha}\}$ for any $i \geq 1$. By \eqref{eq.stepgeneral} and the definition of $s_i$, for all $i\geq 0$,
\begin{equation*}
    \bar{\alpha} = \alpha_{s_{i+1}}  = \alpha_{s_{i}} \tau^{\sum_{k=s_{i}}^{s_{i+1}-1} \Omega_k^g ((1-\Theta_k^g) -\Theta_k^g)} = \bar{\alpha} \tau^{\sum_{k=s_{i}}^{s_{i+1}-1} \Omega_k^g ((1-\Theta_k^g) -\Theta_k^g)},
\end{equation*}
from which it follows that
$0 = \sum_{k=s_{i}}^{s_{i+1}-1} \Omega_k^g ((1-\Theta_k^g) -\Theta_k^g)$.
By the definition of $s_i$, we have: $(i)$ $\alpha_k > \bar{\alpha}$ and $U_k^g = 1$ for any $0 \leq k < s_0$ and $\Omega_k^g = 1$; and, $(ii)$ $\{k: s_i \leq k < s_{i+1}, \Omega_k^g = 1\}$ is subset of either $\{k: s_{i} \leq k < s_{i+1}, \Omega_k^g = 1, U_k^g = 1\}$ or $\{k: s_{i} \leq k < s_{i+1}, \Omega_k^g = 1, U_k^g = 0\}$. As a result, for $i \in \N{}$, 
\begin{equation} \label{eq.bounded.subsummation_be_0}
    \begin{split}
        0 & = \sum_{k=s_{i}}^{s_{i+1}-1} \Omega_k^g ((1-\Theta_k^g) -\Theta_k^g) \\
        & = \sum_{k=s_{i}}^{s_{i+1}-1} \Omega_k^g U_k^g ((1-\Theta_k^g) -\Theta_k^g) = \sum_{k=s_{i}}^{s_{i+1}-1} \Omega_k^g (1-U_k^g) ((1-\Theta_k^g) -\Theta_k^g). 
    \end{split}
\end{equation}
Let $i_t \coloneqq \max\{i \in \N{}: s_i \leq t\}$ for any $t \in \N{}$. We consider two cases.

\emph{Case 1 ($\alpha_t \geq \bar{\alpha}$):} For any  $s_{i_t} \leq k < t$ such that $\Omega_k^g = 1$, it follows that $\alpha_k \geq \bar{\alpha}$ and $U_k^g = 1$. Thus, the exponent in the step size equation \eqref{eq.stepgeneral} can be expressed as
\begin{align}
    \sum_{k=0}^{t-1} \Omega_k^g ((1-\Theta_k^g) -\Theta_k^g) & = \sum_{i=-1}^{i_t} \sum_{k=s_i}^{\min\{s_{i+1}, t\}-1} \Omega_k^g ((1-\Theta_k^g) -\Theta_k^g) \notag \\
    & = \sum_{i=-1}^{i_t} \sum_{k=s_i}^{\min\{s_{i+1}, t\}-1} \Omega_k^g U_k^g ((1-\Theta_k^g) -\Theta_k^g) \label{eq.sum_express1} \\
    & = \sum_{i=0}^{t-1} \Omega_k^g U_k^g ((1-\Theta_k^g) -\Theta_k^g), \notag
\end{align} 
by \eqref{eq.bounded.subsummation_be_0} and the fact that $U_k^g = 1$ for $k=0, \dots, s_0-1$ whenever $\Omega_k^g = 1$. Therefore, $\alpha_t = \alpha_0 \tau^{\sum_{k=0}^{t-1} \Omega_k^g U_k^g((1-\Theta_k^g) -\Theta_k^g)} \geq \bar{\alpha}$. As a result, by \eqref{eq.sum_express1} and the definition of $c_{\tau}$, 
\begin{equation} \label{eq.lg.suc_bound} 
    \sum_{k=0}^{t-1} \Omega_k^g U_k^g (1-\Theta_k^g) \leq \sum_{k=0}^{t-1} \Omega_k^g U_k^g \Theta_k^g + c_{\tau}. 
\end{equation}

\emph{Case 2 ($\alpha_t \leq \bar{\alpha}$):} In this case, for any $s_{i_t} \leq k < t$, it follows that  $\alpha_k \leq \bar{\alpha}$ and $U_k^g = 0$. The step size $\alpha_t$ can be expressed in terms of information at $\alpha_{s_0}$ as follows
\begin{equation} \label{eq.alpha_t_Ug0_form1}
    \alpha_t = \alpha_{s_0}  \tau^{\sum_{k=s_0}^{t-1} \Omega_k^g ((1-\Theta_k^g) -\Theta_k^g)} = \bar{\alpha}  \tau^{\sum_{k=s_0}^{t-1} \Omega_k^g ((1-\Theta_k^g) -\Theta_k^g)}.
\end{equation}
Following the same approach in \eqref{eq.sum_express1}, the exponent in \eqref{eq.alpha_t_Ug0_form1} can be expressed as
\begin{align*}
    & \quad \, \sum_{k=s_0}^{t-1} \Omega_k^g ((1-\Theta_k^g) -\Theta_k^g)  = \sum_{i=0}^{t-1} \Omega_k^g (1-U_k^g) ((1-\Theta_k^g) -\Theta_k^g),
\end{align*}
by \eqref{eq.bounded.subsummation_be_0} and the fact that $1-U_k^g = 0$ for $k=0, \dots, s_0-1$ whenever $\Omega_k^g = 1$. Consequently, $\alpha_t = \bar{\alpha} \cdot \tau^{\sum_{i=0}^{t-1} \Omega_k^g (1-U_k^g) ((1-\Theta_k^g) -\Theta_k^g)} \leq \bar{\alpha}$, and it follows that
\begin{align}\label{eq.sg.suc_bound}
    \sum_{i=0}^{t-1} \Omega_k^g (1-U_k^g) (1-\Theta_k^g) \geq \sum_{i=0}^{t-1} \Omega_k^g (1-U_k^g) \Theta_k^g. 
\end{align}

The bound in \eqref{eq.g.suc_bound_comb} follows by summing \eqref{eq.lg.suc_bound} and \eqref{eq.sg.suc_bound}. The corresponding result for the negative curvature step sizes, \eqref{eq.H.suc_bound_comb}, is obtained analogously.

\eproof

Lemma~\ref{lmm.bounded.bad_upper} bounds the number of iterations corresponding to two ``bad'' scenarios: successful steps with overly small step sizes and unsuccessful steps with overly large step sizes. We now present a lemma characterizing ``good'' iterations.

\blemma \label{lmm.bounded.good_lower}
For all positive integers $t$, and any $0< \bar{p}_g < p_g$ and $0< \bar{p}_H < p_H$ such that $\bar{p}_g  \bar{p}_H + \bar{p}_g + \bar{p}_H -2 > 0$, if $N_{\bar{\epsilon}_g, \bar{\epsilon}_H, \bar{\epsilon}_{\lambda}} > t$ and $\sum_{k=0}^{t-1} I_k^g \geq \bar{p}_g t$,  $\sum_{k=0}^{t-1} I_k^H \geq \bar{p}_H t$, and $\sum_{k=0}^{t-1} I_k^g I_k^H \geq \bar{p}_g \bar{p}_H t$, 
then, 
\begin{equation*}
    \sum_{k=0}^{t-1} I_k^g \Omega_k^g \Theta_k^g U_k^g + \sum_{k=0}^{t-1} I_k^H \Omega_k^H \Theta_k^H U_k^H \geq 
    \tfrac{1}{2} c_{gH}t-c_{\tau},
\end{equation*}
where $c_{gH} \coloneqq \bar{p}_g \bar{p}_H + \bar{p}_g + \bar{p}_H -2$, $c_{\tau} \coloneqq \max\{\log_{\tau} \tfrac{\bar{\alpha}}{\alpha_0}, \log_{\tau} \tfrac{\bar{\beta}}{\beta_0},0\}$ and $\bar\alpha,\bar\beta$ are given in \eqref{def.bar_alpha_beta}. Therefore,
\begin{align}
    
    & \P \left[N_{\bar{\epsilon}_g, \bar{\epsilon}_H, \bar{\epsilon}_{\lambda}} > t, \quad \sum_{k=0}^{t-1} I_k^g \geq \bar{p}_g t, \quad \sum_{k=0}^{t-1} I_k^H \geq \bar{p}_H t, \quad \sum_{k=0}^{t-1} I_k^g I_k^H \geq \bar{p}_g \bar{p}_H t, \right. \label{eq.bounded.good_lower_prob}\\
    & \qquad \qquad \text{and } \ \left. \sum_{k=0}^{t-1} I_k^g \Omega_k^g \Theta_k^g U_k^g + \sum_{k=0}^{t-1} I_k^H \Omega_k^H \Theta_k^H U_k^H < \tfrac{1}{2}c_{gH} t - c_{\tau} \right] = 0. \notag
\end{align} 
\elemma
\bproof
If $N_{\bar{\epsilon}_g, \bar{\epsilon}_H, \bar{\epsilon}_{\lambda}} > t$ and $\sum_{k=0}^{t-1} I_k^g \geq \bar{p}_g t$, $\sum_{k=0}^{t-1} I_k^H \geq \bar{p}_H t$, and $\sum_{k=0}^{t-1} I_k^g I_k^H \geq \bar{p}_g \bar{p}_H t$, then for the descent step, by \eqref{eq.bounded.impossible_case} and \eqref{eq.g.suc_bound_comb},
\begin{align*}
        & \sum_{k=0}^{t-1} I_k^g \Omega_k^g \Theta_k^g U_k^g \\ 
         = &\sum_{k=0}^{t-1} I_k^g \Omega_k^g - \sum_{k=0}^{t-1} I_k^g \Omega_k^g (1-\Theta_k^g) U_k^g - \sum_{k=0}^{t-1} I_k^g \Omega_k^g \Theta_k^g (1-U_k^g) - \sum_{k=0}^{t-1} I_k^g \Omega_k^g (1-\Theta_k^g) (1-U_k^g) \\
        \geq &\sum_{k=0}^{t-1} I_k^g \Omega_k^g - \tfrac{1}{2} \left(\sum_{k=0}^{t-1} \Omega_k^g + c_{\tau} \right) = \tfrac{1}{2} \sum_{k=0}^{t-1} I_k^g \Omega_k^g - \tfrac{1}{2} \sum_{k=0}^{t-1} (1-I_k^g) \Omega_k^g - \tfrac{c_{\tau}}{2}.
\end{align*}
Similarly, for the negative curvature step, by \eqref{eq.bounded.impossible_case} and \eqref{eq.H.suc_bound_comb} it follows that
\begin{equation*}
    \sum_{k=0}^{t-1} I_k^H \Omega_k^H \Theta_k^H U_k^H \geq \tfrac{1}{2} \sum_{k=0}^{t-1} I_k^H \Omega_k^H - \tfrac{1}{2} \sum_{k=0}^{t-1} (1-I_k^H) \Omega_k^H - \tfrac{c_{\tau}}{2}.
\end{equation*}
Combining the above two inequalities, it follows that
\begin{equation} \label{eq.bounded.good_lower}
    \begin{split}
        &\sum_{k=0}^{t-1} I_k^g \Omega_k^g \Theta_k^g U_k^g + \sum_{k=0}^{t-1} I_k^H \Omega_k^H \Theta_k^H U_k^H \\ \geq &\tfrac{1}{2} \sum_{k=0}^{t-1} I_k^g \Omega_k^g + \tfrac{1}{2} \sum_{k=0}^{t-1} I_k^H \Omega_k^H - \tfrac{1}{2} \sum_{k=0}^{t-1} (1-I_k^g) \Omega_k^g - \tfrac{1}{2} \sum_{k=0}^{t-1} (1-I_k^H) \Omega_k^H - c_{\tau}.
    \end{split}
\end{equation}

By $(v)$ of Lemma \ref{lmm.bounded.stocproc} and $I_k^g(1-\Omega_k^g) I_k^H(1-\Omega_k^H) = 0$ for any $k \leq t < N_{\bar{\epsilon}_g, \bar{\epsilon}_H, \bar{\epsilon}_{\lambda}}$, 
\begin{align*}
    \sum_{k=0}^{t-1} I_k^g I_k^H & = \sum_{k=0}^{t-1} I_k^g \Omega_k^g I_k^H \Omega_k^H + \sum_{k=0}^{t-1} I_k^g (1-\Omega_k^g) I_k^H \Omega_k^H + \sum_{k=0}^{t-1} I_k^g \Omega_k^g I_k^H (1-\Omega_k^H) + \sum_{k=0}^{t-1} I_k^g (1-\Omega_k^g) I_k^H (1-\Omega_k^H) \\
    & = \sum_{k=0}^{t-1} I_k^g I_k^H \Omega_k^H + \sum_{k=0}^{t-1} I_k^g \Omega_k^g I_k^H (1-\Omega_k^H) \\
    & \leq \sum_{k=0}^{t-1} I_k^H \Omega_k^H + \sum_{k=0}^{t-1} I_k^g \Omega_k^g.
\end{align*}
Plugging the above into \eqref{eq.bounded.good_lower} and using the definition of $c_{gH}$ completes the proof, 
\begin{align*}
    \sum_{k=0}^{t-1} I_k^g \Omega_k^g \Theta_k^g U_k^g + \sum_{k=0}^{t-1} I_k^H \Omega_k^H \Theta_k^H U_k^H & \geq \tfrac{1}{2} \sum_{k=0}^{t-1} I_k^g I_k^H - \tfrac{1}{2} \sum_{k=0}^{t-1} (1-I_k^g) - \tfrac{1}{2} \sum_{k=0}^{t-1} (1-I_k^H) - c_{\tau} \\
    & \geq \tfrac{1}{2} \bar{p}_g \bar{p}_H t - \tfrac{1}{2}(1-\bar{p}_g) t - \tfrac{1}{2}(1-\bar{p}_H) t -c_{\tau}. 
\end{align*}
\eproof

The above lemma establishes a lower bound, increasing with $t$, on the number of ``good'' iterations (accurate, successful, large step sizes) when the number of accurate iterations is sufficiently large. The main theorem follows.

\btheorem \label{thm.bounded}
     Suppose Assumptions \ref{asm.f.general} and \ref{asm.prob} hold and $e_f \geq 2 \epsilon_f$. Then, for any $0 <\bar{p}_g < p_g$, $0 < \bar{p}_H <p_H$ such that $\bar{p}_g \bar{p}_H + \bar{p}_g + \bar{p}_H -2 \eqqcolon c_{gH} > \tfrac{8e_f}{c_{\bar{\alpha}, \bar{\beta}}}$, 
    and $t \geq T= \tfrac{R}{\tfrac{c_{gH}}{2} - \tfrac{4e_f}{c_{\bar{\alpha}, \bar{\beta}}}}$, 
    \begin{align*}
    \P[N_{\bar{\epsilon}_g, \bar{\epsilon}_H, \bar{\epsilon}_{\lambda}} \leq t] \geq 1- \exp \left( -\tfrac{(p_g - \bar{p}_g)^2}{2p_g^2} t\right) - \exp \left( -\tfrac{(p_H - \bar{p}_H )^2}{2p_H^2} t\right) - \exp \left( -\tfrac{(p_gp_H - \bar{p}_g \bar{p}_H)^2}{2p_g^2p_H^2} t\right),
    \end{align*}
    where $R = \tfrac{f(x_0)-f^*}{c_{\bar{\alpha}, \bar{\beta}}} + c_{\tau}$, $c_{\bar{\alpha}, \bar{\beta}} \coloneqq \min\left\{c_d \bar{\alpha} c_g^2 \bar{\epsilon}_g^2, c_p \bar{\beta}^2 c_H^3 (\max \{\bar{\epsilon}_H, \bar{\epsilon}_{\lambda}\})^3 \right\}$, $\bar\alpha,\bar\beta$ are given in \eqref{def.bar_alpha_beta}, $c_\tau$ is given in Lemma~\ref{lmm.bounded.bad_upper}, 
    $\epsilon_c = 16 \epsilon_f$, and 
    \begin{align}
        \bar{\epsilon}_g & > \max\left\{ \left(\tfrac{\epsilon_c}{(\bar{p}_g \bar{p}_H + \bar{p}_g + \bar{p}_H -2) c_d \bar{\alpha} c_g^2}\right)^{1/2}, \tfrac{2\epsilon_g}{\min\{\eta c_g (1 - \kappa_g), 1-\kappa_g-c_g\}} \right\}, \notag\\
        \bar{\epsilon}_H & > \max\left\{ \left(\tfrac{\epsilon_c}{(\bar{p}_g \bar{p}_H + \bar{p}_g + \bar{p}_H -2) c_p \bar{\beta}^2 c_H^3}\right)^{1/3}, \tfrac{\epsilon_{H}}{c_H} \sqrt{\tfrac{2}{\delta(\gamma- \kappa_H -2 c_p \gamma)}} \right\}, \label{eq.bounded.neibourhood_size}\\
        \bar{\epsilon}_{\lambda} & > \max\left\{ \left(\tfrac{\epsilon_c}{(\bar{p}_g \bar{p}_H + \bar{p}_g + \bar{p}_H -2) c_p \bar{\beta}^2 c_H^3}\right)^{1/3}, \tfrac{\epsilon_{\lambda}}{1-\kappa_{\lambda}-c_H} \right\}. \notag
    \end{align}
\etheorem

\bproof
    When $N_{\bar{\epsilon}_g, \bar{\epsilon}_H, \bar{\epsilon}_{\lambda}} > t \geq T$, by Lemma \ref{lmm.bounded.good_upper} and the definition of $T$ 
\begin{align*}
    \sum_{k=0}^{t-1} I_k^g \Omega_k^g \Theta_k^g U_k^g + \sum_{k=0}^{t-1} I_k^H \Omega_k^H \Theta_k^H U_k^H  \leq \tfrac{f(x_0)-f(x_t)}{c_{\bar{\alpha}, \bar{\beta}}} + \tfrac{4e_f}{c_{\bar{\alpha}, \bar{\beta}}} t & < \tfrac{f(x_0)-f^*}{c_{\bar{\alpha}, \bar{\beta}}} + \tfrac{1}{2} c_{gH} t - R \\
    & = \tfrac{1}{2} c_{gH} t - c_{\tau}.
\end{align*}
We denote the events
\begin{align*}
    E_g \coloneqq \left\{\sum_{k=0}^{t-1} I_k^g \geq \bar{p}_g t \right\}, \quad E_H \coloneqq \left\{\sum_{k=0}^{t-1} I_k^H \geq \bar{p}_H t \right\},  \quad E_{gH} \coloneqq \left\{\sum_{k=0}^{t-1} I_k^g I_k^H \geq \bar{p}_g \bar{p}_H t \right\}.
\end{align*}
The probability of not reaching a second-order stationary point by iteration~$t$ is
\begin{align*}
    &\P[N_{\bar{\epsilon}_g, \bar{\epsilon}_H, \bar{\epsilon}_{\lambda}} > t] \\
    = &\P \left[N_{\bar{\epsilon}_g, \bar{\epsilon}_H, \bar{\epsilon}_{\lambda}} > t,\ \sum_{k=0}^{t-1} I_k^g \Omega_k^g \Theta_k^g U_k^g + \sum_{k=0}^{t-1} I_k^H \Omega_k^H \Theta_k^H U_k^H < \tfrac{1}{2} c_{gH} t - c_{\tau} \right] \\
    = &\P \left[N_{\bar{\epsilon}_g, \bar{\epsilon}_H, \bar{\epsilon}_{\lambda}} > t,\ \sum_{k=0}^{t-1} I_k^g \Omega_k^g \Theta_k^g U_k^g + \sum_{k=0}^{t-1} I_k^H \Omega_k^H \Theta_k^H U_k^H < \tfrac{1}{2} c_{gH} t - c_{\tau}, \ E_g \cap E_H \cap E_{gH} \right] \\
    & \; + \P \left[N_{\bar{\epsilon}_g, \bar{\epsilon}_H, \bar{\epsilon}_{\lambda}} > t,\ \sum_{k=0}^{t-1} I_k^g \Omega_k^g \Theta_k^g U_k^g + \sum_{k=0}^{t-1} I_k^H \Omega_k^H \Theta_k^H U_k^H < \tfrac{1}{2} c_{gH} t - c_{\tau}, \ E_g^c \cup E_H^c \cup E_{gH}^c \right] \\
    \leq &0 + \P \left\{ E_g^c \cup E_H^c \cup E_{gH}^c \right\} \\
    \leq &\P[E_g^c] + \P[E_H^c] + \P[E_{gH}^c] \\
    \leq & \exp \left( -\tfrac{(p_g - \bar{p}_g)^2}{2p_g^2} t\right) + \exp \left( -\tfrac{(p_H - \bar{p}_H)^2}{2p_H^2} t\right) + \exp \left( -\tfrac{(p_g p_H - \bar{p}_g \bar{p}_H)^2}{2p_g^2p_H^2} t\right),
\end{align*}
where the first inequality follows by \eqref{eq.bounded.good_lower_prob} and the last inequality by a special case of Lemma \ref{lmm.azuma} where $p_f = \bar{p}_f = 1$ and $I_k^f = \hat{I}_k^f = 1$ for all positive integers $k$.
\eproof

\bremark
We make the following remarks about Theorem~\ref{thm.bounded}. 
\begin{itemize}[leftmargin=0.5cm]
    \item The dependence of the complexity bounds on $\bar{\epsilon}_g$ and $\bar{\epsilon}_H$ parallels the deterministic case \cite{berahas2024exploiting, curtis2019exploiting, li2025randomized}, where iterates converge to a neighborhood of a stationary point. The neighborhood parameters $\bar{\epsilon}_g$, $\bar{\epsilon}_H$, and $\bar{\epsilon}_{\lambda}$ depend on $\epsilon_f$ with different orders, satisfying $\bar{\epsilon}_g > \mathcal{O}(\epsilon_f^{1/2} + \epsilon_g)$, $\bar{\epsilon}_H > \mathcal{O}(\epsilon_f^{1/3} + \epsilon_H)$, and $\bar{\epsilon}_{\lambda} > \mathcal{O}(\epsilon_f^{1/3} + \epsilon_{\lambda})$. Consequently, the achievable second-order accuracy of Algorithm~\ref{alg.theory} is $\left(\mathcal{O}(\epsilon_f^{1/2} + \epsilon_g),\, 
    \mathcal{O}(\epsilon_f^{1/3} + \epsilon_H),\, 
    \mathcal{O}(\epsilon_f^{1/3} + \epsilon_{\lambda})\right)$.
    \item When $p_g = p_H = 1$ and $\epsilon_f = \epsilon_g = \epsilon_H = \epsilon_{\lambda} = 0$, the result reduces to the deterministic inexact setting \cite[Theorem 2.10]{berahas2024exploiting}, where the lower bounds on $\bar{\epsilon}_g$ and $\bar{\epsilon}_H$ vanish. Moreover, if $p_H = 1$, the requirement on $p_g$ becomes $p_g > 1/2$, matching the first-order condition in \cite{jin2021high}.
    \item The probability that $N_{\bar{\epsilon}_g, \bar{\epsilon}_H, \bar{\epsilon}_{\lambda}}$ exceeds $t$ decays exponentially in $t$ for $t \ge T$, where $T = \mathcal{O}\!\left(\epsilon_f^{-1} + \epsilon_g^{-2} + (\max\{\epsilon_H, \epsilon_{\lambda}\})^{-3}\right)$.
\end{itemize}
\eremark

\subsection{Convergence analysis: Subexponential noise case} \label{subsec.analysis.subexp}

We now extend the analysis to the case where the noise in the function estimates has a subexponential tail, corresponding to Oracle \ref{def.0th-order_orc}.\ref{orc.0th-order_subexp_orc}. The lemmas from the bounded-noise setting still hold after replacing $I_k^g (p_g)$ with $I_k^f I_k^g (p_f p_g)$ and $I_k^H (p_H)$ with $\hat{I}_k^f I_k^H (p_f p_H)$. We first state the analogue of Lemma \ref{lmm.bounded.stocproc} for this setting.
\blemma \label{lmm.subexp.stocproc}
    Let $e_f \geq 4/a + 2\epsilon_f$ and the neighborhood parameters ($\bar{\epsilon}_g, \bar{\epsilon}_H, \bar{\epsilon}_{\lambda}$) set as prescribed in \eqref{eq.value.eps}.
    Then, there exist constants $\bar{\alpha}, \bar{\beta} > 0$  \eqref{def.bar_alpha_beta} and nondecreasing functions $h_d(\cdot), h_p(\cdot): \mathbb{R} \rightarrow \mathbb{R}$, satisfying $h_d(\alpha), h_p(\beta) > 0$ for any $\alpha, \beta > 0$, such that for any realization of Algorithm \ref{alg.theory} the following hold for all $k < N_{\bar{\epsilon}_g, \bar{\epsilon}_H, \bar{\epsilon}_{\lambda}}$:
    \begin{enumerate}[leftmargin=0.5cm]
        \item [(i)] If the function and gradient estimates at the descent step of iteration $k$ are accurate (i.e., $I_k^f I_k^g = 1$) and $\|g_k\|_2 \geq c_g \bar{\epsilon}_g$ (i.e., $\Omega_k^g = 1$), then for  $\alpha_k \leq \bar{\alpha}$, the iteration (descent step) is successful ($\Theta_k^g = 1$), which implies $\alpha_{k+1} = \tau^{-1} \alpha_k$.
        \item [(ii)] If the function and Hessian estimates at the negative curvature step of iteration $k$ are accurate (i.e., $\hat{I}_k^f I_k^H = 1$) and $\lambda_k \leq -c_H \max\{\bar{\epsilon}_H, \bar{\epsilon}_{\lambda}\}$ (i.e., $\Omega_k^H = 1$), then for $\beta_k \leq \bar{\beta}$, the iteration (negative curvature step) is successful ($\Theta_k^H = 1$), which implies $\beta_{k+1} = \tau^{-1} \beta_k$.
        \item [(iii)] When $I_k^fI_k^g \Omega_k^g \Theta_k^g = 1$, $f(x_{k+1}) \leq f(x_k) - h_d(\alpha_k) + 3e_f + \hat{e}_k + \max\{\hat{e}_k^+, \hat{e}_k^-\}$.
        \item [(iv)] When $\hat{I}_k^f I_k^H \Omega_k^H \Theta_k^H = 1$, $f(x_{k+1}) \leq f(x_k) - h_p(\beta_k) + 3e_f + e_k + e_k^+$.
        \item [(v)] $I_k^g(1-\Omega_k^g) I_k^H(1-\Omega_k^H) = 0$.
    \end{enumerate}
    In summary,
    \begin{align} \label{eq.subexp.summary_bounded}
        f(x_{k+1}) \leq \begin{cases}
            f(x_k) - h_d(\alpha_k) + 3e_f + \hat{e}_k + \max\{\hat{e}_k^+, \hat{e}_k^-\}, & \text{if } I_k^f I_k^g \Omega_k^g \Theta_k^g = 1, \\
            f(x_k) - h_p(\beta_k) + 3e_f + e_k + e_k^+, & \text{if } \hat{I}_k^f I_k^H \Omega_k^H \Theta_k^H = 1, \\
            f(x_k) + 2e_f+ e_k + e_k^+ + \hat{e}_k + \max\{\hat{e}_k^+, \hat{e}_k^-\}, & \text{otherwise}.
        \end{cases}
    \end{align}
\elemma
\bproof
The proof follows the same arguments as Lemma~\ref{lmm.bounded.stocproc}. In particular, part~$(v)$ is identical in both statement and proof and is therefore omitted. We restate the remaining arguments below, highlighting the differences.

For the descent step, when $I_k^f I_k^g \Omega_k^g = 1$, the case becomes the same as in $(i)$ of Lemma \ref{lmm.bounded.stocproc}, and by \eqref{eq.bounded.fk+bound} and \eqref{eq.bounded.Fk+bound}, as long as $\alpha_k \leq \bar{\alpha}$, defined in \eqref{def.bar_alpha_beta}, it follows that
\begin{align*}
        F(x_k + \alpha_k d_k, \xi_k^{(0+)}) & \leq F(x_k, \xi_k^{(0)}) + c_d \alpha_k g_k^{\mathsf{T}} d_k + e_k + e_k^+ \leq F(x_k, \xi_k^{(0)}) + c_d \alpha_k g_k^{\mathsf{T}} d_k + e_f,
\end{align*}
which implies that the descent step is successful ($\Theta_k^g = 1$). When $I_k^f I_k^g \Omega_k^g \Theta_k^g = 1$, the case is similar to $(iii)$ of Lemma \ref{lmm.bounded.stocproc}, and
\begin{equation*}
    \begin{split}
        f(x_k + \alpha_k d_k) & \leq f(x_k) -h_d(\alpha_k) + e_f + e_k + e_k^+ \leq f(x_k) -h_d(\alpha_k) + 2e_f, 
    \end{split}
\end{equation*}
where $h_d(\alpha) \coloneqq c_g^2 \bar{\epsilon}_g^2 \alpha$. The worst-case bound for $f(\hat{x}_k)$, analogous to \eqref{eq.bounded.dc_bound_bad}, is
\begin{align*}
    f(\hat{x}_k) = f(x_k + \alpha_k d_k) \leq f(x_k) + e_f + e_k + e_k^+.
\end{align*}

For the negative curvature step, when $\hat{I}_k^f I_k^H \Omega_k^H = 1$, for $\beta_k \leq \bar{\beta}$ defined in \eqref{def.bar_alpha_beta}, 
\begin{align*}
    &\min\{F(\hat{x}_k + \beta_k q_k, \hat{\xi}_k^{(0+)}), F(\hat{x}_k - \beta_k q_k, \hat{\xi}_k^{(0-)}) \} \\ 
    \leq &\min\{f(\hat{x}_k + \beta_k q_k), f(\hat{x}_k - \beta_k q_k) \} + \max\{\hat{e}_k^+, \hat{e}_k^-\} \\
    \leq &f(\hat{x}_k) + c_p \beta_k^2 q_k^{\mathsf{T}} H_k q_k + \max\{\hat{e}_k^+, \hat{e}_k^-\} \\
     \leq & F(\hat{x}_k, \hat{\xi}_k^{(0)}) + c_p \beta_k^2 q_k^{\mathsf{T}} H_k q_k + \hat{e}_k + \max\{\hat{e}_k^+, \hat{e}_k^-\} \leq F(\hat{x}_k, \hat{\xi}_k^{(0)}) + c_p \beta_k^2 q_k^{\mathsf{T}} H_k q_k + e_f
\end{align*}
and the step is successful ($\Theta_k^H = 1$). When $\hat{I}_k^f I_k^H \Omega_k^H\Theta_k^H = 1$, 
\begin{equation*}
    \begin{split}
        f(\hat{x}_k + \beta_k p_k) & \leq \min\{F(\hat{x}_k + \beta_k q_k, \hat{\xi}_k^{(0+)}), F(\hat{x}_k - \beta_k q_k, \hat{\xi}_k^{(0-)}) \} + \max \{\hat{e}_k^+, \hat{e}_k^-\} \\
        & \leq F(\hat{x}_k, \hat{\xi}_k^{(0)}) + c_p \beta_k^2 q_k^{\mathsf{T}} H_k q_k + e_f + \max \{\hat{e}_k^+, \hat{e}_k^-\} \\
        & \leq f(\hat{x}_k) + \hat{e}_k - h_p(\beta_k) + e_f + \max \{\hat{e}_k^+, \hat{e}_k^-\} \leq f(\hat{x}_k) - h_p(\beta_k) + 2e_f,
    \end{split}
\end{equation*}
where $h_p(\beta) \coloneqq c_p \gamma \delta^2 c_H^3 (\max \{\bar{\epsilon}_H, \bar{\epsilon}_{\lambda}\})^3 \beta^2$. The bound, analogous to  \eqref{eq.bounded.nc_bound_bad}, is
\begin{align*}
    f(x_{k+1}) = f(\hat{x}_k + \beta_k p_k) \leq f(\hat{x}_k) + e_f + \hat{e}_k + \max \{\hat{e}_k^+, \hat{e}_k^-\}.
\end{align*}
Combining the inequalities for different cases above, yields the desired results.
\eproof 

The following lemma, analogous to Lemma~\ref{lmm.bounded.good_upper}, bounds the number of “good” iterations and follows directly from Lemma~\ref{lmm.subexp.stocproc}.

\blemma \label{lmm.subexp.good_upper}
For any positive integer $t$, we have
\begin{align*}
    \sum_{k=0}^{t-1} I_k^f I_k^g \Omega_k^g \Theta_k^g U_k^g + \sum_{k=0}^{t-1} \hat{I}_k^f I_k^H \Omega_k^H \Theta_k^H U_k^H & \leq \tfrac{f(x_0) - f(x_t)}{c_{\bar{\alpha}, \bar{\beta}}} + \tfrac{3e_f}{c_{\bar{\alpha}, \bar{\beta}}} t \\
    & \quad + \tfrac{1}{c_{\bar{\alpha}, \bar{\beta}}} \sum_{k=0}^{t-1} \left(e_k + e_k^+ + \hat{e}_k + \max\{\hat{e}_k^+, \hat{e}_k^-\} \right),
\end{align*}
where $c_{\bar{\alpha}, \bar{\beta}} \coloneqq \min\{h_d(\bar{\alpha}), h_p(\bar{\beta})\}$ and $\bar\alpha,\bar\beta$ are given in \eqref{def.bar_alpha_beta}.
\elemma
\bproof Similar to the proof of Lemma~\ref{lmm.bounded.good_upper}, by incorporating the step size indicator variables $U_k^g$ and $U_k^H$, one can re-write \eqref{eq.subexp.summary_bounded} as 
$$
    f(x_{k+1}) \leq \begin{cases}
       f(x_k) - h_d(\bar{\alpha}) + 3e_f + \hat{e}_k + \max\{\hat{e}_k^+, \hat{e}_k^-\}, & \text{if } I_k^f I_k^g \Omega_k^g \Theta_k^g U_k^g = 1, \\
       f(x_k) - h_p(\bar{\beta}) + 3e_f + e_k + e_k^+, & \text{if } \hat{I}_k^f I_k^H \Omega_k^H \Theta_k^H U_k^H = 1, \\
       f(x_k) + 3e_f + e_k + e_k^+ + \hat{e}_k + \max\{\hat{e}_k^+, \hat{e}_k^-\}, & \text{otherwise}.
    \end{cases}
$$
since $h_d(\cdot)$ and $h_p(\cdot)$ are non-decreasing functions on $\mathbb{R}_{\geq 0}$.
Summing the inequalities above from $k=0$ to $t-1$, re-arranging the terms, and using the definition of $c_{\bar{\alpha}, \bar{\beta}}$ completes the proof.
\eproof

Lemma~\ref{lmm.bounded.bad_upper} applies directly in this setting. The next lemma, analogous to Lemma~\ref{lmm.bounded.good_lower}, provides a lower bound on the number of ``good'' iterations (accurate, successful, large step sizes), implying the proportion remains bounded away from zero as $t$ increases.

\blemma \label{lmm.subexp.good_lower}
For all positive integers $t$, and any $0 < \bar{p}_{fg} < p_f p_g$ and $0 < \bar{p}_{fH} < p_f p_H$ such that $\bar{p}_{fg} \cdot \bar{p}_{fH} + \bar{p}_{fg} + \bar{p}_{fH} -2 > 0$, if $N_{\bar{\epsilon}_g, \bar{\epsilon}_H, \bar{\epsilon}_{\lambda}} > t$ and $\sum_{k=0}^{t-1} I_k^f I_k^g \geq \bar{p}_g t$, $\sum_{k=0}^{t-1} \hat{I}_k^f I_k^H \geq \bar{p}_{fH} t$, and $\sum_{k=0}^{t-1} I_k^f I_k^g \hat{I}_k^f I_k^H \geq \bar{p}_{fg} \bar{p}_{fH} t$, then, 
\begin{equation*}
    \sum_{k=0}^{t-1} I_k^f I_k^g \Omega_k^g \Theta_k^g U_k^g + \sum_{k=0}^{t-1} \hat{I}_k^f I_k^H \Omega_k^H \Theta_k^H U_k^H \geq \tfrac{1}{2} (\bar{p}_{fg} \bar{p}_{fH} + \bar{p}_{fg} + \bar{p}_{fH} -2) t - c_{\tau} = \tfrac{1}{2} c_{gH}t-c_{\tau},
\end{equation*}
where $c_{gH} \coloneqq \bar{p}_{fg} \bar{p}_{fH} + \bar{p}_{fg} + \bar{p}_{fH} -2$, $c_{\tau} \coloneqq \max\{\log_{\tau} \tfrac{\bar{\alpha}}{\alpha_0}, \log_{\tau} \tfrac{\bar{\beta}}{\beta_0},0\}$ and $\bar\alpha,\bar\beta$ are given in \eqref{def.bar_alpha_beta}. 
Therefore,
\begin{equation*}
    \begin{split}
        & \P \left[N_{\bar{\epsilon}_g, \bar{\epsilon}_H, \bar{\epsilon}_{\lambda}} > t, \quad \sum_{k=0}^{t-1} I_k^f I_k^g \geq \bar{p}_{fg} t, \quad \sum_{k=0}^{t-1} \hat{I}_k^f I_k^H \geq \bar{p}_{fH} t, \quad \sum_{k=0}^{t-1} I_k^f I_k^g \hat{I}_k^f I_k^H \geq \bar{p}_{fg} \bar{p}_{fH} t, \right. \\
        & \qquad \qquad \text{and } \ \left. \sum_{k=0}^{t-1} I_k^f I_k^g \Omega_k^g \Theta_k^g U_k^g + \sum_{k=0}^{t-1} \hat{I}_k^f I_k^H \Omega_k^H \Theta_k^H U_k^H < \tfrac{1}{2}c_{gH} t - c_{\tau} \right] = 0.
    \end{split}
\end{equation*}
\elemma
\bproof
The proof parallels that of Lemma~\ref{lmm.bounded.good_lower}, incorporating $I_k^f$ and $\hat{I}_k^f$ and their associated probabilities, and is therefore omitted.
\eproof

The main result follows and uses arguments similar to those of Theorem~\ref{thm.bounded}.
\btheorem \label{thm.prob}
Suppose Assumptions \ref{asm.f.general} and \ref{asm.prob} hold and $e_f \geq 2\epsilon_f + 5/a$. Then, for any $s \geq  0$, $0 < \bar{p}_{fg} < p_f p_g$, $0 < \bar{p}_{fH} < p_f p_H$ such that
$\bar{p}_{fg} \bar{p}_{fH} + \bar{p}_{fg} + \bar{p}_{fH} -2 \eqqcolon c_{gH} > \tfrac{10e_f + 4s}{c_{\bar{\alpha}, \bar{\beta}}}$, and $t \geq T \coloneqq \tfrac{R}{\tfrac{c_{gH}}{2} - \tfrac{5e_f + 2s}{c_{\bar{\alpha}, \bar{\beta}}}}$, 
\begin{align*}
    \P\{N_{\bar{\epsilon}_g, \bar{\epsilon}_H, \bar{\epsilon}_{\lambda}} \leq t\}  \geq 1 &- \exp \left( -\tfrac{(p_f p_g - \bar{p}_{fg})^2}{2p_f^2p_g^2} t\right) - \exp \left( -\tfrac{(p_f p_H - \bar{p}_{fH} )^2}{2p_f^2 p_H^2} t\right) \\
     \quad \ \ &- \exp \left( -\tfrac{(p_f^2 p_gp_H - \bar{p}_{fg} \bar{p}_{fH})^2}{2p_f^4p_g^2p_H^2} t\right) - 2 \exp\left(-\tfrac{a}{4}s t \right)
\end{align*}
where $R = \tfrac{f(x_0)-f^*}{c_{\bar{\alpha}, \bar{\beta}}} + c_{\tau}$, $c_{\bar{\alpha}, \bar{\beta}} \coloneqq \min\left\{c_d \bar{\alpha} c_g^2 \bar{\epsilon}_g^2, c_p \bar{\beta}^2 c_H^3 \left(\max\{\bar{\epsilon}_H, \bar{\epsilon}_{\lambda}\}\right)^ 3 \right\}$, $\bar\alpha,\bar\beta$ are given in \eqref{def.bar_alpha_beta}, $c_\tau$ is given in Lemma~\ref{lmm.bounded.bad_upper}, and the bounds for $\bar{\epsilon}_g$, $\bar{\epsilon}_H$, and $\bar{\epsilon}_{\lambda}$ are in the same form as in \eqref{eq.bounded.neibourhood_size} with factor $\epsilon_c = 16 \epsilon_f + 32/a + 4s$.
\etheorem
\bproof
We define the following two events: $A_t \coloneqq \left\{ \tfrac{1}{t} \sum_{k=0}^{t-1} (e_k + e_k^+) \leq e_f + s \right\}$ and $\hat{A}_t \coloneqq \left\{ \tfrac{1}{t} \sum_{k=0}^{t-1} \left(\hat{e}_k + \max\{\hat{e}_k^+, \hat{e}_k^-\} \right) \leq e_f + s \right\}$. 
The event $\{N_{\bar{\epsilon}_g, \bar{\epsilon}_H, \bar{\epsilon}_{\lambda}} > t\}$ can then be divided into two events using $A_t$ and $\hat{A}_t$,
\begin{equation} \label{eq.subexp.prob.overt_2events}
    \P \left[N_{\bar{\epsilon}_g, \bar{\epsilon}_H, \bar{\epsilon}_{\lambda}} > t \right] = \P \left[N_{\bar{\epsilon}_g, \bar{\epsilon}_H, \bar{\epsilon}_{\lambda}} > t, \ A_t \cap \hat{A}_t \right] + \P \left[ N_{\bar{\epsilon}_g, \bar{\epsilon}_H, \bar{\epsilon}_{\lambda}} > t, \ A_t^c \cup \hat{A}_t^c \right].
\end{equation}
We first bound the latter term. By $e_f \geq 2\epsilon_f + 5/a$ and Lemma \ref{lmm.prob.ef_lg}, it follows that
\begin{align} \label{eq.subexp.prob_bound1}
        & \quad \, \, \P \left[ N_{\bar{\epsilon}_g, \bar{\epsilon}_H, \bar{\epsilon}_{\lambda}} > t, \ A_t^c \cup \hat{A}_t^c \right] \\
        & \leq \P\left[A_t^c \right] + \P \left[\hat{A}_t^c \right] \nonumber\\
        & = \P\left[\tfrac{1}{t} \sum_{k=0}^{t-1} (e_k + e_k^+) > e_f + s \right] + \P\left[\tfrac{1}{t} \sum_{k=0}^{t-1} (\hat{e}_k + \max\{\hat{e}_k^+, \hat{e}_k^-\}) > e_f + s \right] \nonumber\\
        & \leq \P\left[\tfrac{1}{t} \sum_{k=0}^{t-1} (e_k + e_k^+) > 2\epsilon_f + \tfrac{5}{a} + s \right] + \P\left[\tfrac{1}{t} \sum_{k=0}^{t-1} (\hat{e}_k + \max\{\hat{e}_k^+, \hat{e}_k^-\}) > 2\epsilon_f + \tfrac{5}{a} + s \right] \nonumber\\
        & \leq 2 \exp \left( -\tfrac{a}{4} st \right) \nonumber
\end{align}
The former term in \eqref{eq.subexp.prob.overt_2events} can be decomposed as follows, 
\begin{equation}
\begin{aligned}  \label{eq.subexp.prob.bounded_2events}
    & \quad \, \, \P \left[N_{\bar{\epsilon}_g, \bar{\epsilon}_H, \bar{\epsilon}_{\lambda}} > t, \ A_t \cap \hat{A}_t \right] \\
    &  = \P \left[N_{\bar{\epsilon}_g, \bar{\epsilon}_H, \bar{\epsilon}_{\lambda}} > t, \ A_t \cap \hat{A}_t, \sum_{k=0}^{t-1} I_k^f I_k^g \Omega_k^g \Theta_k^g U_k^g + \sum_{k=0}^{t-1} \hat{I}_k^f I_k^H \Omega_k^H \Theta_k^H U_k^H < \tfrac{1}{2}c_{gH} t - c_{\tau} \right] \\
    & \quad + \P \left[N_{\bar{\epsilon}_g, \bar{\epsilon}_H, \bar{\epsilon}_{\lambda}} > t, \ A_t \cap \hat{A}_t, \sum_{k=0}^{t-1} I_k^f I_k^g \Omega_k^g \Theta_k^g U_k^g + \sum_{k=0}^{t-1} \hat{I}_k^f I_k^H \Omega_k^H \Theta_k^H U_k^H \geq \tfrac{1}{2}c_{gH} t - c_{\tau} \right] 
\end{aligned}
\end{equation}
By Lemma \ref{lmm.subexp.good_upper}, when $A_t \cap \hat{A}_t$ is true and $N_{\bar{\epsilon}_g, \bar{\epsilon}_H, \bar{\epsilon}_{\lambda}} > t \geq T$, it follows that
\begin{align*} 
    \sum_{k=0}^{t-1} I_k^f I_k^g \Omega_k^g \Theta_k^g U_k^g + \sum_{k=0}^{t-1} \hat{I}_k^f I_k^H \Omega_k^H \Theta_k^H U_k^H & \leq \tfrac{f(x_0) - f(x_t)}{c_{\bar{\alpha}, \bar{\beta}}} + \tfrac{3e_f}{c_{\bar{\alpha}, \bar{\beta}}} t \\
    & \quad + \tfrac{1}{c_{\bar{\alpha}, \bar{\beta}}} \sum_{k=0}^{t-1} \left(e_k + e_k^+ + \hat{e}_k + \max\{\hat{e}_k^+, \hat{e}_k^-\} \right) \\
    & \leq \tfrac{f(x_0)-f^*}{c_{\bar{\alpha}, \bar{\beta}}} + \tfrac{5e_f + 2s}{c_{\bar{\alpha}, \bar{\beta}}} t  < \tfrac{c_{gH}}{2}t-c_{\tau} .
\end{align*}
Thus, the second term in \eqref{eq.subexp.prob.bounded_2events} is 0, and we only need to bound the first term,
\begin{equation} \label{eq.subexp.prob.simp_bound}
    \begin{split}
        & \P \left[N_{\bar{\epsilon}_g, \bar{\epsilon}_H, \bar{\epsilon}_{\lambda}} > t, \ A_t \cap \hat{A}_t, \ \sum_{k=0}^{t-1} I_k^f I_k^g \Omega_k^g \Theta_k^g U_k^g + \sum_{k=0}^{t-1} \hat{I}_k^f I_k^H \Omega_k^H \Theta_k^H U_k^H < \tfrac{1}{2}c_{gH} t - c_{\tau} \right] \\
        \leq \: & \P \left[N_{\bar{\epsilon}_g, \bar{\epsilon}_H, \bar{\epsilon}_{\lambda}} > t, \ \sum_{k=0}^{t-1} I_k^f I_k^g \Omega_k^g \Theta_k^g U_k^g + \sum_{k=0}^{t-1} \hat{I}_k^f I_k^H \Omega_k^H \Theta_k^H U_k^H < \tfrac{1}{2}c_{gH} t - c_{\tau} \right].
    \end{split}
\end{equation}
The rest of the proof is similar to the one in the bounded-noise setting. We denote the events
$$
    E_{fg} = \left\{\sum_{k=0}^{t-1} I_k^f I_k^g \geq \bar{p}_{fg} t \right\}, \ \ E_{fH} = \left\{\sum_{k=0}^{t-1} \hat{I}_k^f I_k^H \geq \bar{p}_{fH} t \right\}, \ \ E_{fgH} = \left\{\sum_{k=0}^{t-1} I_k^f I_k^g \hat{I}_k^f I_k^H \geq \bar{p}_{fg} \bar{p}_{fH} t \right\}.
$$
Then the probability in \eqref{eq.subexp.prob.simp_bound} can be bounded as follows 
\begin{align} \label{eq.subexp.prob_bound2}
        & \quad \ \P \left[N_{\bar{\epsilon}_g, \bar{\epsilon}_H, \bar{\epsilon}_{\lambda}} > t, \ \sum_{k=0}^{t-1} I_k^f I_k^g \Omega_k^g \Theta_k^g U_k^g + \sum_{k=0}^{t-1} \hat{I}_k^f I_k^H \Omega_k^H \Theta_k^H U_k^H < \tfrac{1}{2}c_{gH} t - c_{\tau} \right] \\
        & = \P \left[N_{\bar{\epsilon}_g, \bar{\epsilon}_H, \bar{\epsilon}_{\lambda}} > t,\ \sum_{k=0}^{t-1} I_k^f I_k^g \Omega_k^g \Theta_k^g U_k^g + \sum_{k=0}^{t-1} \hat{I}_k^f I_k^H \Omega_k^H \Theta_k^H U_k^H < \tfrac{1}{2} c_{gH} t - c_{\tau}, \ E_{fg} \cap E_{fH} \cap E_{fgH} \right] \nonumber\\
        & \quad + \P \left[N_{\bar{\epsilon}_g, \bar{\epsilon}_H, \bar{\epsilon}_{\lambda}} > t,\ \sum_{k=0}^{t-1} I_k^f I_k^g \Omega_k^g \Theta_k^g U_k^g + \sum_{k=0}^{t-1} \hat{I}_k^f I_k^H \Omega_k^H \Theta_k^H U_k^H < \tfrac{1}{2} c_{gH} t - c_{\tau}, \ E_{fg}^c \cup E_{fH}^c \cup E_{fgH}^c \right] \nonumber\\
        & = 0 + \P \left[N_{\bar{\epsilon}_g, \bar{\epsilon}_H, \bar{\epsilon}_{\lambda}} > t,\ \sum_{k=0}^{t-1} I_k^f I_k^g \Omega_k^g \Theta_k^g U_k^g + \sum_{k=0}^{t-1} \hat{I}_k^f I_k^H \Omega_k^H \Theta_k^H U_k^H < \tfrac{1}{2} c_{gH} t - c_{\tau}, \ E_{fg}^c \cup E_{fH}^c \cup E_{fgH}^c \right] \nonumber\\
        & \leq \P \left[ E_{fg}^c \cup E_{fH}^c \cup E_{fgH}^c \right] \nonumber\\
        & \leq \P\left[E_{fg}^c\right] + \P\left[E_{fH}^c\right] + \P\left[E_{fgH}^c\right] \nonumber\\
        & \leq \exp \left( -\tfrac{(p_f p_g - \bar{p}_{fg})^2}{2p_f^2 p_g^2} t\right) + \exp \left( -\tfrac{(p_f p_H - \bar{p}_{fH})^2}{2p_f^2p_H^2} t\right) + \exp \left( -\tfrac{(p_f^2 p_g p_H - \bar{p}_{fg} \bar{p}_{fH})^2}{2p_f^4 p_g^2p_H^2} t\right). \nonumber
\end{align}
where the fourth line follows by Lemma \ref{lmm.subexp.good_lower} and the last inequality follows by Lemma \ref{lmm.azuma}.
Combining \eqref{eq.subexp.prob_bound1} and \eqref{eq.subexp.prob_bound2}, it follows that
\begin{align*}
    \P\left[ N_{\bar{\epsilon}_g, \bar{\epsilon}_H, \bar{\epsilon}_{\lambda}} > t\right] & \leq 2 \exp \left( - \tfrac{a}{4} st \right) + \exp \left( -\tfrac{(p_f p_g - \bar{p}_{fg})^2}{2p_f^2 p_g^2} t\right) + \exp \left( -\tfrac{(p_f p_H - \bar{p}_{fH})^2}{2p_f^2p_H^2} t\right) \\
    & \quad + \exp \left( -\tfrac{(p_f^2 p_g p_H - \bar{p}_{fg} \bar{p}_{fH})^2}{2p_f^4 p_g^2p_H^2} t\right).
\end{align*}
Re-arranging the above completes the proof.
\eproof
\bremark
The theorem parallels Theorem~\ref{thm.bounded} for bounded noise, with two modifications reflecting the subexponential tails of the function oracle:
\begin{itemize}[leftmargin=0.5cm]
    \item The noise-control parameter $e_f$ becomes $2\epsilon_f + 5/a > 2\epsilon_f$, where $a>0$ and $\epsilon_f \ge 0$ are oracle parameters. The lower bounds on the neighborhood parameters $\bar{\epsilon}_g$, $\bar{\epsilon}_H$, and $\bar{\epsilon}_{\lambda}$ now depend on $e_f$ and a free parameter $s \ge 0$, satisfying $
        \bar{\epsilon}_g > \mathcal{O} ((\epsilon_f + a^{-1}+s)^{1/2} + \epsilon_g), \ \bar{\epsilon}_H > \mathcal{O} ((\epsilon_f + a^{-1}+s)^{1/3}, \epsilon_H), \ \bar{\epsilon}_{\lambda} > \mathcal{O} ((\epsilon_f + a^{-1}+s)^{1/3} + \epsilon_{\lambda})
    $. Consequently, the achievable second-order accruacy of Algorithm~\ref{alg.theory} is $\left(\mathcal{O}((\epsilon_f + a^{-1}+s)^{1/2} + \epsilon_g), \mathcal{O}((\epsilon_f + a^{-1}+s)^{1/3} + \epsilon_H), \mathcal{O}((\epsilon_f + a^{-1}+s)^{1/3} + \epsilon_{\lambda})\right)$.
    \item In addition to the three exponential terms present in the bounded-noise analysis (now involving $p_f p_g$ and $p_f p_H$ due to the function oracle in the line search), an additional factor $2\exp(-\tfrac{a}{4} s t)$ appears, capturing the subexponential tail of the function noise. The choice of $s>0$ trades a slightly larger neighborhood for faster decay of the failure probability. As before, the probability that $N_{\bar{\epsilon}_g, \bar{\epsilon}_H, \bar{\epsilon}_{\lambda}}$ exceeds $t$ decays exponentially in $t$ for $t \ge T$, where $T = \mathcal{O}\!\left(\epsilon_f^{-1} + a + s^{-1} + \epsilon_g^{-2} + (\max\{\epsilon_H, \epsilon_{\lambda}\})^{-3}\right)$. 
\end{itemize}
\eremark

\section{Numerical Experiments} \label{sec.experiment}

\begin{figure}[H]
    \centering
    \begin{subfigure}{1\textwidth}
        \includegraphics[width=0.24\linewidth]{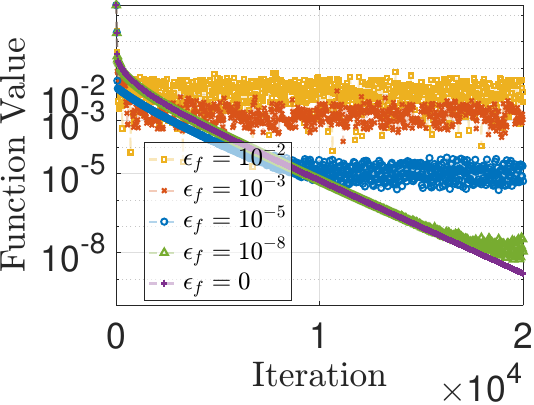}
        \includegraphics[width=0.24\linewidth]{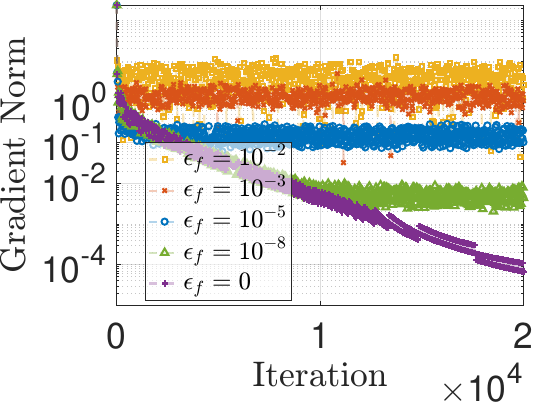}
        \includegraphics[width=0.24\linewidth]{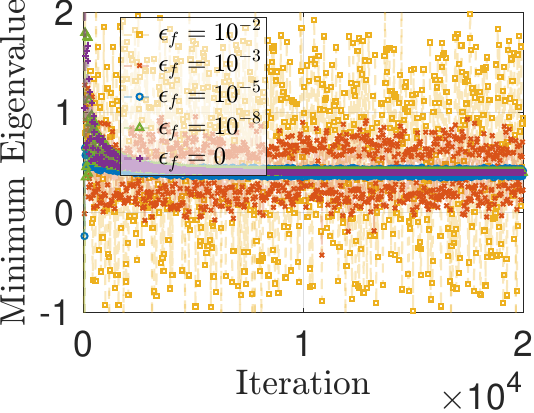}
        \includegraphics[width=0.24\linewidth]{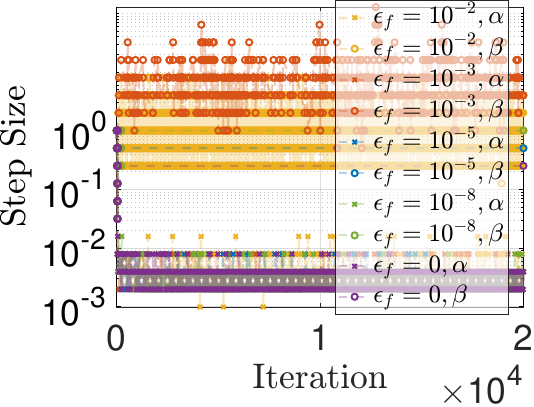}

        \includegraphics[width=0.24\linewidth]{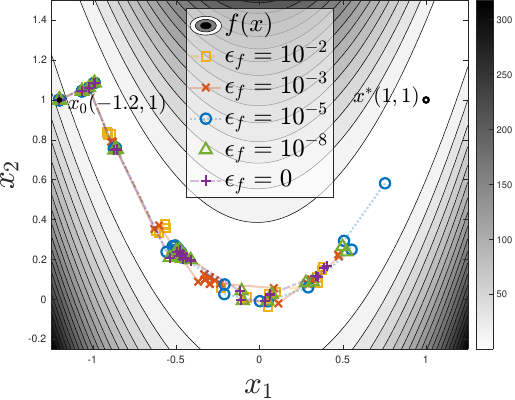}
        \includegraphics[width=0.24\linewidth]{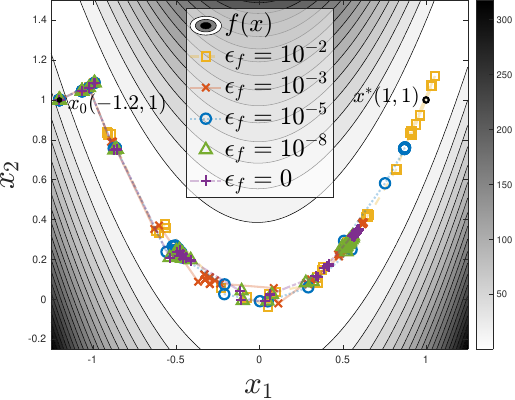}
        \includegraphics[width=0.24\linewidth]{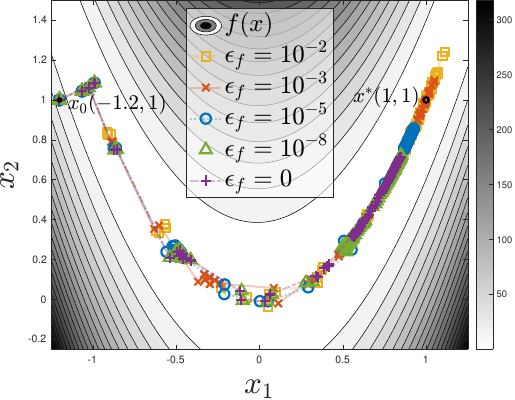}
        \includegraphics[width=0.24\linewidth]{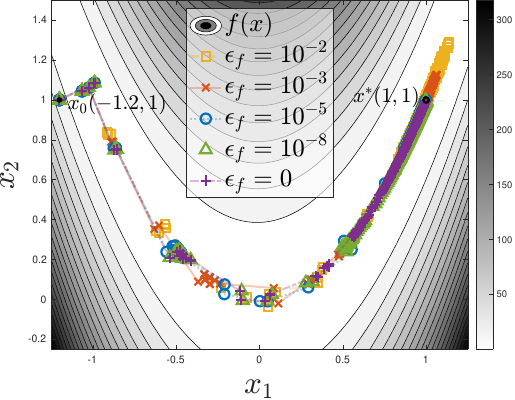}
    \end{subfigure}
    \caption{{\small Sensitivity of Algorithm~\ref{alg.theory} on the Rosenbrock problem for $\epsilon_f \in \{10^{-2}, 10^{-3}, 10^{-5}, 10^{-8}, 0\}$ with $e_f = 2\epsilon_f$. Top row: metrics vs.\ iterations. Bottom row: contour plots with trajectories (after 100, 200, 2000, and 20000 iterations).}}
    \label{fig.ROSENBR.sensitivity.eps_f}
\end{figure}

\begin{figure}[H]
    \centering
    \begin{subfigure}{1\textwidth}
        \includegraphics[width=0.25\linewidth]{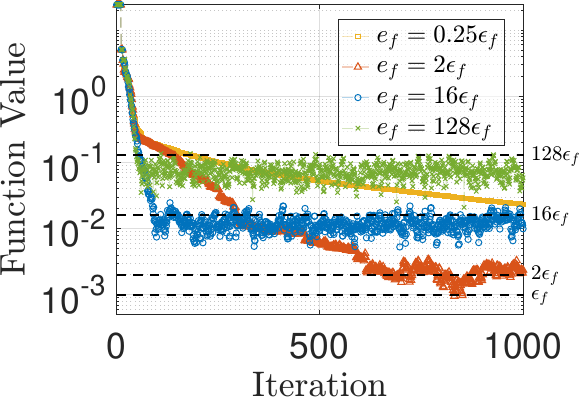}
        \includegraphics[width=0.245\linewidth]{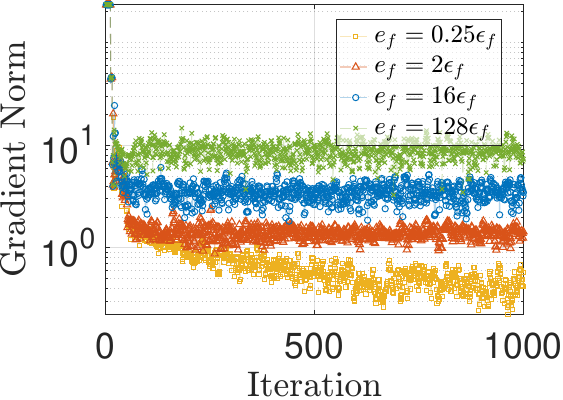}
        \hspace{-0.5em}
        \includegraphics[width=0.25\linewidth]{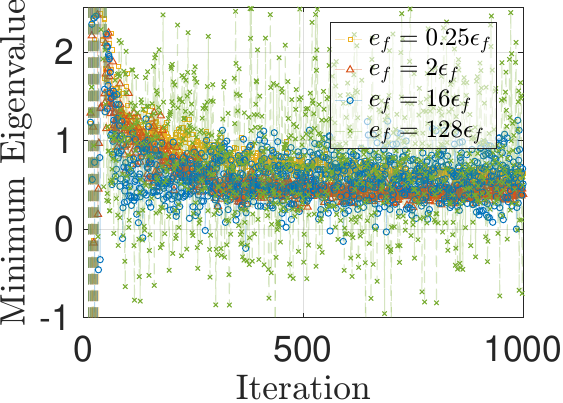}
        \hspace{-0.5em}
        \includegraphics[width=0.25\linewidth]{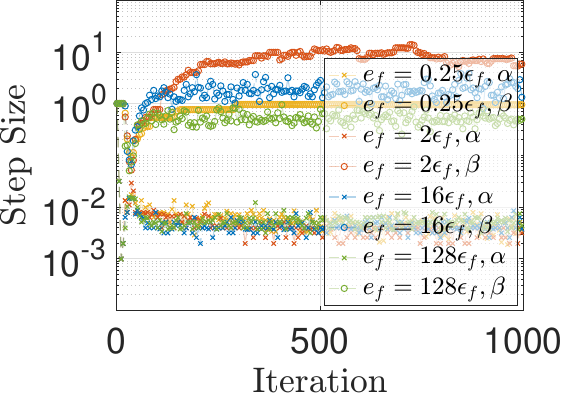} \\ \vspace{-0.2cm}

        \includegraphics[width=0.24\linewidth]{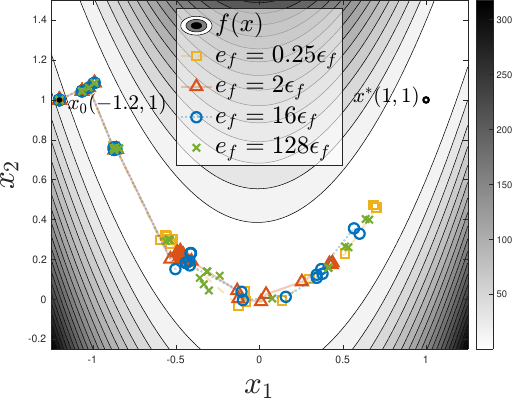}
        \includegraphics[width=0.24\linewidth]{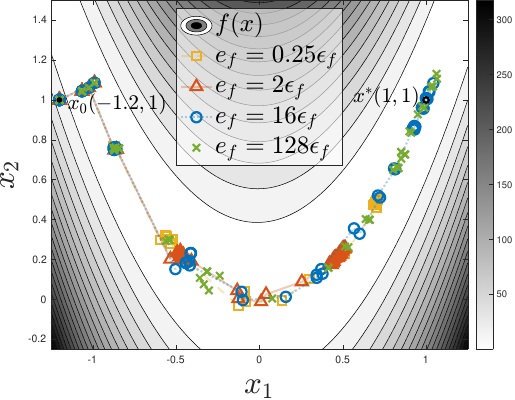}
        \includegraphics[width=0.24\linewidth]{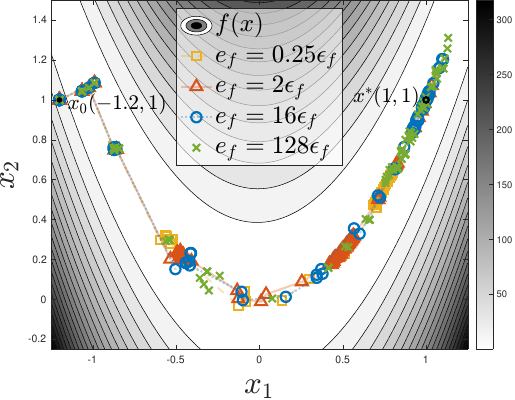}
        \includegraphics[width=0.24\linewidth]{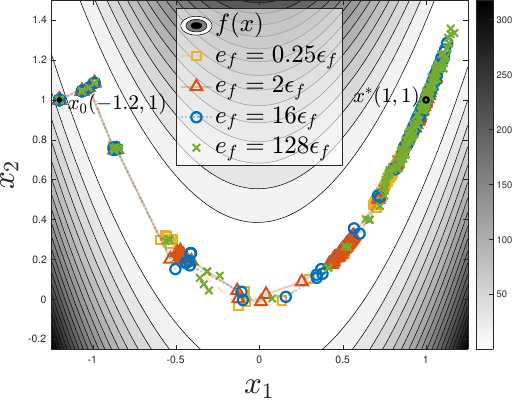}
    \end{subfigure}
    \caption{{\small 
    Sensitivity of Algorithm~\ref{alg.theory} on the Rosenbrock problem for $\epsilon_f = 10^{-3}$ and $e_f \in \{0.25\epsilon_f, 2\epsilon_f, 16\epsilon_f, 128\epsilon_f\}$. Results averaged over 10 runs. Top row: metrics vs. iterations. Bottom row: contour plots with trajectories (after 100, 200, 1000, and 5000 iterations).
    }}
    \label{fig.ROSENBR.sensitivity.e_f}
\end{figure}

\begin{figure}[H]
    \centering
    \begin{subfigure}{1\textwidth}
        \includegraphics[width=0.25\linewidth]{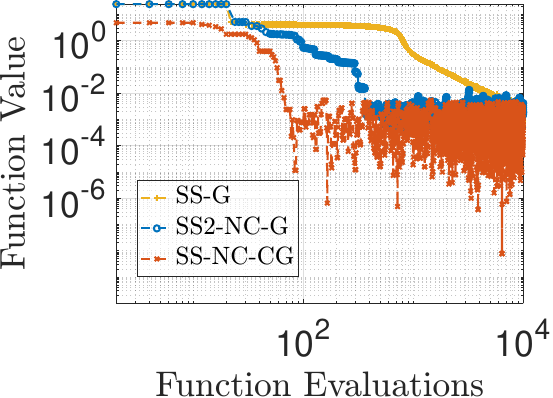}
        \hspace{-0.5em}
        \includegraphics[width=0.25\linewidth]{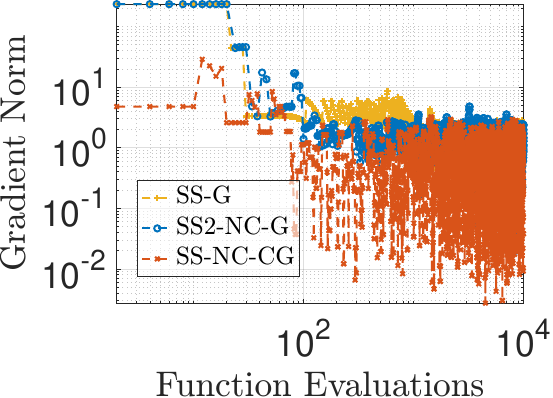}
        \hspace{-0.5em}
        \includegraphics[width=0.25\linewidth]{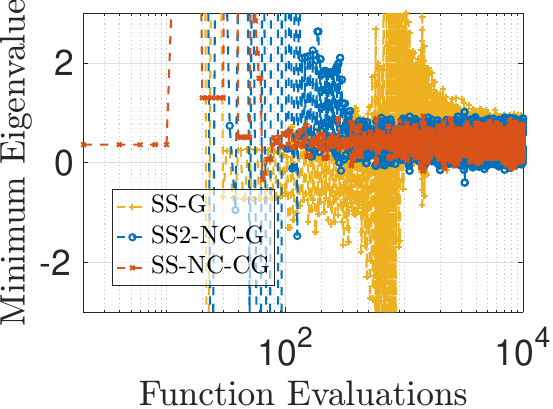}
        \hspace{-0.5em}
        \includegraphics[width=0.25\linewidth]{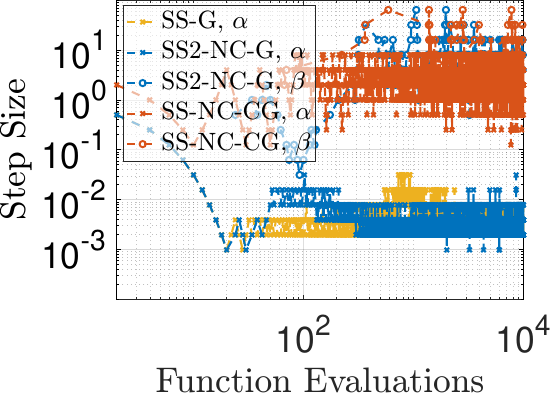} \\

    \includegraphics[width=0.24\linewidth]{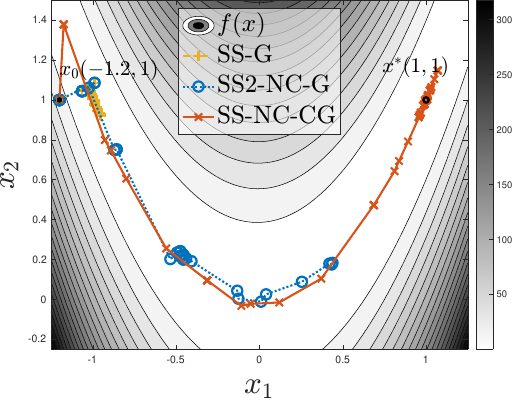}
        \includegraphics[width=0.24\linewidth]{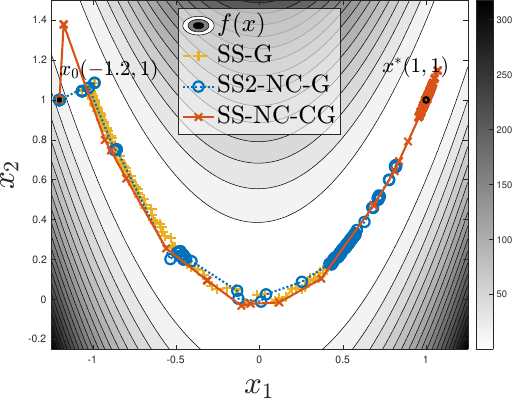}
        \includegraphics[width=0.24\linewidth]{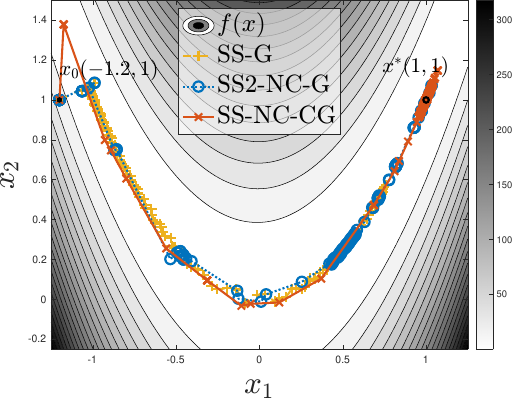}
        \includegraphics[width=0.24\linewidth]{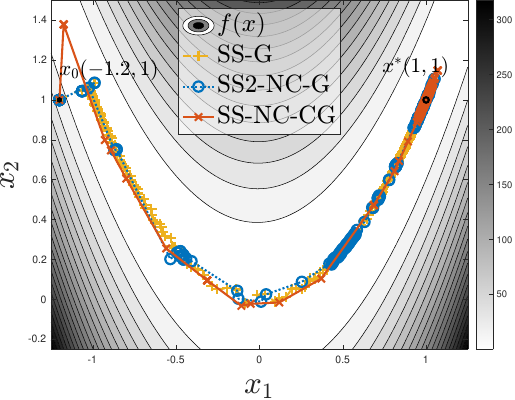}
    \end{subfigure}
    \caption{{\small 
    Comparison of Algorithm~\ref{alg.theory} (\texttt{SS2-NC-G}), \texttt{SS-G}, and \texttt{SS-NC-CG} on the Rosenbrock problem with $\epsilon_f = 10^{-3}$ and $e_f = 2\epsilon_f$. Top row: metrics vs. function evaluations. Bottom row: contour plots with iterate trajectories (after 100, 500, 1000, and 5000 iterations).
    }}
    \label{fig.ROSENBR.1e-3.2}
\end{figure}

In this section, we report numerical experiments evaluating the practical performance of the proposed two-step method (Algorithm~\ref{alg.theory}, denoted \texttt{SS2-NC-G}). First, we examine the sensitivity of the method to key parameters, namely the function-noise level $\epsilon_f$ (Figure~\ref{fig.ROSENBR.sensitivity.eps_f}) and the Armijo noise-tolerance parameter $e_f$ (Figure~\ref{fig.ROSENBR.sensitivity.e_f}), and then compare with Adaptive Line-search with Oracle Estimations (denoted \texttt{SS-G})~\cite{jin2021high} and a conjugate-gradient-based variant (denoted \texttt{SS-NC-CG})~\cite{royer2020newton} (Figure~\ref{fig.ROSENBR.1e-3.2}). We use the Rosenbrock function~\cite{gratton2025s2mpj} as a test case, and evaluate performance in terms of objective function value, gradient norm, minimum Hessian eigenvalue, and step size with respect to iterations and function evaluations.

To simulate oracle outputs, we perturb the exact function, gradient, and Hessian values with controlled bounded noise. This corresponds to Oracle~\ref{def.0th-order_orc}.\ref{orc.0th-order_bounded}, Oracle~\ref{def.1st-order_orc} with $p_g = 1$ and $\kappa_g = 0$, and Oracle~\ref{def.2nd-order_orc} with $p_H = 1$ and $\kappa_H = \kappa_{\lambda} = 0$. To simplify the experimental setting and align with the theoretical scaling, we impose $\epsilon_f = \epsilon_g^2 = \epsilon_H^3$. The resulting perturbation models are described below.
\begin{itemize}[leftmargin=0.5cm]
    \item \textbf{Function estimates:} Noisy function evaluations are generated by adding bounded random noise $\epsilon_f$ to the exact objective value, $F \leftarrow f(x) + \epsilon_f  \mathcal{U}(-1,1)$. 
    \item \textbf{Gradient estimates:} To generate gradient perturbations with bounded norm, we first draw $r \sim \mathcal{N}(0,I_n)$ and normalize it as $u \leftarrow r/\|r\|_2$, yielding a random direction uniformly distributed on the unit sphere. We then draw $U \sim \mathcal{U}(0,1)$ and set the radius $\rho \leftarrow \epsilon_g  U^{1/n}$. The noisy gradient estimate is $g \leftarrow \nabla f(x) + \rho  u$. This construction guarantees that the perturbation satisfies $\|g - \nabla f(x)\|_2 \leq \epsilon_g$.
    \item \textbf{Hessian estimates:} Similarly, we generate Hessian perturbations by first drawing a random matrix $R \sim \mathcal{N}(0,1)^{n \times n}$ and normalizing it as $u \leftarrow R/\|R\|_2$. We then draw $U \sim \mathrm{Uniform}(0,1)$ and set $\rho \leftarrow \epsilon_H  U^{1/n^2}$. The noisy Hessian estimate is given by $H \leftarrow \nabla^2 f(x) + \rho u$. Under this construction, the perturbation satisfies $\|H - \nabla^2 f(x)\|_2 \leq \epsilon_H$.
\end{itemize}

For all methods, we use the same step-search hyperparameters $\alpha_0 = \beta_0 = 1$, $\tau = 0.5$, $c_d = c_p = 0.2$, and set $c_g = 0$. For negative curvature methods, we set the early-termination threshold to 
$c_H \max \{\bar{\epsilon}_H, \bar{\epsilon}_{\lambda}\} = 10^{-3}$. 
All experiments were conducted in MATLAB~R2021b. 

We first investigate the sensitivity of \texttt{SS2-NC-G} to the noise level $\epsilon_f \in \{10^{-2}, 10^{-3}, 10^{-5}, 10^{-8}, 0\}$. As shown in Figure~\ref{fig.ROSENBR.sensitivity.eps_f}, the noise magnitude determines the neighborhood of convergence. Smaller values of $\epsilon_f$ (e.g., $0$ and $10^{-8}$) produce smoother trajectories and convergence to smaller neighborhoods, albeit with slower initial progress, whereas larger values yield more variability and larger neighborhoods, often with faster initial decrease. These trends are consistent with the theoretical scaling $\epsilon_f = \epsilon_g^2 = \epsilon_H^3$, under which the iteration complexity before entering the high-probability regime is $\mathcal{O}(\epsilon_f^{-1})$. Although larger-noise runs may progress more rapidly at early stages, smaller-noise runs ultimately attain better objective values and stationarity. The contour plots further illustrate the differing trajectory behavior across noise levels.

We next examine the effect of misspecifying $e_f$, which governs step acceptance under noisy function evaluations. We test
$e_f \in \{0.25\epsilon_f, 2\epsilon_f, 16\epsilon_f, 128\epsilon_f\}$.
In the bounded-noise setting (Oracle~\ref{def.0th-order_orc}.\ref{orc.0th-order_bounded}), the theory requires $e_f \ge 2\epsilon_f$ for convergence; nevertheless, we also consider smaller values to evaluate behavior in the regime in which the noise is underestimated. As shown in Figure~\ref{fig.ROSENBR.sensitivity.e_f}, larger $e_f$ yields convergence to a larger neighborhood, consistent with the dependence of the accuracy bounds on $e_f$, while permitting larger step sizes and faster initial progress. The choice $e_f = 2\epsilon_f$, motivated by the theory, seems to balance early progress and final accuracy.

In the final experiment, we compare \texttt{SS2-NC-G} with the first-order method \texttt{SS-G} and \texttt{SS-NC-CG}. The latter follows a simplified version of the framework in \cite{royer2020newton} and computes either a Newton-type step or a negative curvature step, with step sizes determined by step search due to the stochastic nature of the problem. We report the evolution of the objective value, gradient norm, minimum eigenvalue, and step size with respect to iterations and function evaluations. Results on the Rosenbrock problem with $\epsilon_f = 10^{-3}$ and $e_f = 2\epsilon_f$ are shown in Figure~\ref{fig.ROSENBR.1e-3.2}. Methods incorporating negative curvature (\texttt{SS2-NC-G} and \texttt{SS-NC-CG}) reduce the objective more effectively in regions of negative curvature, while \texttt{SS-G} exhibits slower progress near saddle-like regions. The trajectory plots reflect this behavior.

\section{Final Remarks} \label{sec.remarks}

We consider optimization problems in which the objective function and its derivatives are corrupted by noise and accessed through probabilistic oracles. Within this framework, we have developed a two-step negative curvature method incorporating a step-search procedure with a relaxed Armijo-type sufficient-decrease condition and a novel mechanism for selecting the negative curvature direction. Under mild assumptions on the accuracy of the probabilistic oracles, the proposed approach is endowed with high-probability second-order convergence and iteration-complexity guarantees. Numerical experiments on a standard nonconvex problem illustrate the efficiency and robustness of the proposed framework.

\newpage
\bibliographystyle{plain}
\bibliography{references}

\begin{thebibliography}{10}

\bibitem{allen2018natasha}
Zeyuan Allen-Zhu.
\newblock Natasha 2: Faster non-convex optimization than sgd.
\newblock {\em Advances in neural information processing systems}, 31, 2018.

\bibitem{azuma1967weighted}
Kazuoki Azuma.
\newblock Weighted sums of certain dependent random variables.
\newblock {\em Tohoku Mathematical Journal, Second Series}, 19(3):357--367,
  1967.

\bibitem{bellavia2022linesearch}
Stefania Bellavia, Eugenio Fabrizi, and Benedetta Morini.
\newblock Linesearch newton-cg methods for convex optimization with noise.
\newblock {\em Annali dell'Universita'di Ferrara}, 68(2):483--504, 2022.

\bibitem{berahas2024exploiting}
Albert~S Berahas, Raghu Bollapragada, and Wanping Dong.
\newblock Exploiting negative curvature in conjunction with adaptive sampling:
  theoretical results and a practical algorithm.
\newblock {\em Computational Optimization and Applications}, pages 1--37, 2026.

\bibitem{berahas2019derivative}
Albert~S Berahas, Richard~H Byrd, and Jorge Nocedal.
\newblock Derivative-free optimization of noisy functions via quasi-newton
  methods.
\newblock {\em SIAM Journal on Optimization}, 29(2):965--993, 2019.

\bibitem{berahas2021global}
Albert~S Berahas, Liyuan Cao, and Katya Scheinberg.
\newblock Global convergence rate analysis of a generic line search algorithm
  with noise.
\newblock {\em SIAM Journal on Optimization}, 31(2):1489--1518, 2021.

\bibitem{berahas2025sequential}
Albert~S Berahas, Miaolan Xie, and Baoyu Zhou.
\newblock A sequential quadratic programming method with high-probability
  complexity bounds for nonlinear equality-constrained stochastic optimization.
\newblock {\em SIAM Journal on Optimization}, 35(1):240--269, 2025.

\bibitem{bollapragada2018adaptive}
Raghu Bollapragada, Richard Byrd, and Jorge Nocedal.
\newblock Adaptive sampling strategies for stochastic optimization.
\newblock {\em SIAM Journal on Optimization}, 28(4):3312--3343, 2018.

\bibitem{bottou2018optimization}
L{\'e}on Bottou, Frank~E Curtis, and Jorge Nocedal.
\newblock Optimization methods for large-scale machine learning.
\newblock {\em SIAM review}, 60(2):223--311, 2018.

\bibitem{byrd2012sample}
Richard~H Byrd, Gillian~M Chin, Jorge Nocedal, and Yuchen Wu.
\newblock Sample size selection in optimization methods for machine learning.
\newblock {\em Mathematical programming}, 134(1):127--155, 2012.

\bibitem{cao2023first}
Liyuan Cao, Albert~S Berahas, and Katya Scheinberg.
\newblock First-and second-order high probability complexity bounds for
  trust-region methods with noisy oracles.
\newblock {\em Math. Programming}, pages 1--52, 2023.

\bibitem{carter1991global}
Richard~G Carter.
\newblock On the global convergence of trust region algorithms using inexact
  gradient information.
\newblock {\em SIAM Journal on Numerical Analysis}, 28(1):251--265, 1991.

\bibitem{cartis2011adaptive}
Coralia Cartis, Nicholas~IM Gould, and Philippe~L Toint.
\newblock Adaptive cubic regularisation methods for unconstrained optimization.
  {{Part I}}: Motivation, convergence and numerical results.
\newblock {\em Mathematical Programming}, 127(2):245--295, 2011.

\bibitem{cartis2018global}
Coralia Cartis and Katya Scheinberg.
\newblock Global convergence rate analysis of unconstrained optimization
  methods based on probabilistic models.
\newblock {\em Mathematical Programming}, 169:337--375, 2018.

\bibitem{curtis2019exploiting}
Frank~E Curtis and Daniel~P Robinson.
\newblock Exploiting negative curvature in deterministic and stochastic
  optimization.
\newblock {\em Mathematical Programming}, 176(1):69--94, 2019.

\bibitem{forsgren1995computing}
Anders Forsgren, Philip~E Gill, and Walter Murray.
\newblock Computing modified newton directions using a partial cholesky
  factorization.
\newblock {\em SIAM Journal on Scientific Comp.}, 16(1):139--150, 1995.

\bibitem{friedlander2012hybrid}
Michael~P Friedlander and Mark Schmidt.
\newblock Hybrid deterministic-stochastic methods for data fitting.
\newblock {\em SIAM Journal on Scientific Computing}, 34(3):A1380--A1405, 2012.

\bibitem{goldfarb1980curvilinear}
Donald Goldfarb.
\newblock Curvilinear path steplength algorithms for minimization which use
  directions of negative curvature.
\newblock {\em Mathematical programming}, 18(1):31--40, 1980.

\bibitem{gratton2025s2mpj}
Serge Gratton and Ph~L Toint.
\newblock S2mpj and cutest optimization problems for matlab, python and julia.
\newblock {\em Optimization Methods and Software}, pages 1--33, 2025.

\bibitem{jin2021high}
Billy Jin, Katya Scheinberg, and Miaolan Xie.
\newblock High probability complexity bounds for line search based on
  stochastic oracles.
\newblock {\em Neural Information Processing Systems}, 34:9193--9203, 2021.

\bibitem{li2025randomized}
Shuyao Li and Stephen~J Wright.
\newblock A randomized algorithm for nonconvex minimization with inexact
  evaluations and complexity guarantees.
\newblock {\em JOTA}, 207(3):66, 2025.

\bibitem{mccormick1977modification}
Garth~P McCormick.
\newblock A modification of armijo's step-size rule for negative curvature.
\newblock {\em Mathematical programming}, 13(1):111--115, 1977.

\bibitem{more1979use}
Jorge~J Mor{\'e} and Danny~C Sorensen.
\newblock On the use of directions of negative curvature in a modified newton
  method.
\newblock {\em Mathematical Programming}, 16:1--20, 1979.

\bibitem{nesterov2006cubic}
Yurii Nesterov and Boris~T Polyak.
\newblock Cubic regularization of newton method and its global performance.
\newblock {\em Mathematical programming}, 108(1):177--205, 2006.

\bibitem{nocedal1999numerical}
Jorge Nocedal and Stephen~J Wright.
\newblock {\em Numerical optimization}.
\newblock Springer, 1999.

\bibitem{paquette2020stochastic}
Courtney Paquette and Katya Scheinberg.
\newblock A stochastic line search method with expected complexity analysis.
\newblock {\em SIAM Journal on Optimization}, 30(1):349--376, 2020.

\bibitem{pasupathy2018sampling}
Raghu Pasupathy, Peter Glynn, Soumyadip Ghosh, and Fatemeh~S Hashemi.
\newblock On sampling rates in simulation-based recursions.
\newblock {\em SIAM Journal on Optimization}, 28(1):45--73, 2018.

\bibitem{rockafellar2000optimization}
R~Tyrrell Rockafellar, Stanislav Uryasev, et~al.
\newblock Optimization of conditional value-at-risk.
\newblock {\em Journal of risk}, 2:21--42, 2000.

\bibitem{royer2020newton}
Cl{\'e}ment~W Royer, Michael O’Neill, and Stephen~J Wright.
\newblock A newton-cg algorithm with complexity guarantees for smooth
  unconstrained optimization.
\newblock {\em Math. Programming}, 180:451--488, 2020.

\bibitem{royer2018complexity}
Cl{\'e}ment~W Royer and Stephen~J Wright.
\newblock Complexity analysis of second-order line-search algorithms for smooth
  nonconvex optimization.
\newblock {\em SIAM Journal on Optimization}, 28(2):1448--1477, 2018.

\bibitem{vershynin2018high}
Roman Vershynin.
\newblock {\em High-dimensional probability: An introduction with applications
  in data science}, volume~47.
\newblock Cambridge university press, 2018.

\bibitem{xu2020newton}
Peng Xu, Fred Roosta, and Michael~W Mahoney.
\newblock {Newton-type methods for non-convex optimization under inexact
  Hessian information}.
\newblock {\em Mathematical Programming}, 184(1):35--70, 2020.

\end{thebibliography}


\end{document}